\documentclass[12pt]{article}
\usepackage{amsmath, amssymb, amsthm, amscd}
\usepackage{graphics}

\newtheorem{theorem}      {Theorem}
\newtheorem*{theorem*}    {Theorem}

\newtheorem{lemma}        {Lemma}

\newtheorem{proposition}  {Proposition}
\newtheorem{corollary}    {Corollary}
\numberwithin{corollary}{theorem}

\renewcommand{\qed}{$\blacksquare$}

\newtheorem*{claim}       {Claim}

\theoremstyle{definition}
\newtheorem{example}      {Example}
\newtheorem{ex}[example]      {Example}
\newtheorem{remark}       {Remark}
\newtheorem{definition}   {Definition}

\newcommand{\vol}{\mathrm{Vol}}
\newcommand{\cvx}{\mathrm{Conv}}
\newcommand{\diag}{\mathrm{Diag}}

\newcommand{\supp}{\mathrm{Supp}}
\newcommand{\linear}{\mathrm{Lin}}
\newcommand{\punct}[1] {\hspace{0.3cm}\text{#1}}
\newcommand{\toric}[1] {\ensuremath{\mathcal T^{#1}}}
\newcommand{\defeq}{\stackrel{\scriptstyle\mathrm{def}}{=}}
\newcommand{\romone}{\scriptscriptstyle\mathrm{I}}
\newcommand{\romtwo}{\scriptscriptstyle\mathrm{I\!I}}
\newcommand{\romthree}{\scriptscriptstyle\mathrm{I\!I\!I}}

\newcommand{\mymu}{\ensuremath{\boldsymbol{\mu}}}

\newcommand{\finsler}[1]{\ensuremath{
\left| \! \left| \! \left| #1 \right| \! \right| \! \right|}}

\newcommand{\thth}{\mathrm{\underline{th}}} 
\newcommand{\eps}{\varepsilon} 
\newcommand{\cF}{\mathcal{F}} 
\newcommand{\Pro}{\mathbb{P}} 
\newcommand{\Z}{\mathbb{Z}} 
\newcommand{\R}{\mathbb{R}} 
\newcommand{\C}{\mathbb{C}} 
\newcommand{\Cn}{\C^n} 

\newcounter{formula}
\newlength{\formulaparindent}
\numberwithin{formula}{subsection}

\newlength{\formulawidth}
\newcommand{\formula}[2]{
\ 
\refstepcounter{formula}
\vspace{1ex} \\
\fbox{\parbox{\formulawidth}{\vspace{2ex}\hspace{\formulaparindent} Formula \theformula: #1
\[ {#2} \] }}
}

\numberwithin{equation}{subsection}

\newcount\n
\newcount\m
\def\twodigits#1{\ifnum #1<10 0\fi \number#1}

\def\hours{\n=\time \divide\n 60
  \m=-\n \multiply\m 60 \advance\m \time
    \twodigits\n :\twodigits\m}

\author{Gregorio Malajovich\thanks{Departamento de Matem\'atica Aplicada,
Universidade Federal do Rio de Janeiro, Caixa Postal 68530, CEP 21945--970,
Rio de Janeiro, RJ, Brasil. {\tt http://www.labma.ufrj.br/\~{}gregorio} ~. 
{\bf e-mail:} {\tt gregorio@labma.ufrj.br} ~.
On leave at the Department of Mathematics, City University of Hong Kong, 
83 Tat Chee Ave, Kowloon, Hong Kong.}
\thanks{G.M.'s visit to City University of Hong Kong was supported
by CERG grants ~9040393--730, 9040402--730, and~9040188 }  
\and J. Maurice Rojas\thanks{Department of Mathematics, City University
of Hong Kong (83 Tat Chee Ave., Kowloon, HONG KONG) and 
Department of Mathematics, Texas A\&M University 
(College Station, Texas 77843-3368, USA). {\tt 
http://math.cityu.edu.hk/~mamrojas} ~. {\bf e-mail:} 
{\tt mamrojas@cityu.edu.hk} (before January 2001), {\tt 
rojas@math.tamu.edu} (after January 2001). }
\,\thanks{Partially supported by Hong Kong UGC grant \#9040402-730 
and a US National Science Foundation Mathematical Sciences 
Postdoctoral Fellowship.}
}

\title{
\mbox{}\\
\vspace{-4cm}
Random Sparse Polynomial Systems  
} 

\date{December 11, 2000}

\begin{document}
\enlargethispage{1cm}

\maketitle

\mbox{}\hspace{3cm}{\sc {\small To Steve Smale on his 70$^\thth$ birthday.}} 

\thispagestyle{empty}
\begin{abstract}
  Let $f\!:=\!(f^1,\ldots,f^n)$ be a sparse random polynomial system. This 
means that each $f^i$ has fixed support (list of possibly non-zero 
coefficients) and each coefficient has a Gaussian probability 
distribution of arbitrary variance.
\par
  We express the expected number of roots of $f$ inside a region $U$ as 
the integral over $U$ of a certain {\bf mixed volume} 
form. When $U = (\mathbb C^*)^n$, the classical mixed volume is recovered.
\par
  The main result in this paper is a bound on the probability that the 
condition number of $f$
on the region $U$ is larger than $1/\eps$. This bound depends on
the integral of the mixed volume form over $U$, and on a certain 
intrinsic invariant of $U$ as a subset of a toric manifold.
\par
  Polynomials with real coefficients are also considered, and bounds
for the expected number of real roots and for the condition number
are given. 
\par
  The connection between zeros of sparse random polynomial systems,
K\"ahler geometry, and mechanics (momentum maps) is discussed. 
\end{abstract}


\noindent 
{\bf 2000 Math Subject Classification:} 
65H10, 
52A39. 

\newpage
\vspace{\stretch{1}} \ \\
\tableofcontents
\vspace{\stretch{1}}
\newpage

\section{Introduction}

A complexity theory of homotopy algorithms for solving dense systems
of polynomial equations was developed in ~\cite{BEZ1,BEZ2,BEZ3,BEZ4,BEZ5}.
(see also~\cite[Ch.~8--14]{BCSS}). The number of steps for these homotopy 
algorithms was bounded in terms of a condition number and the
B\'ezout number. 

One of the main features of that theory was unitary invariance: the roots
of a dense system of polynomial equations are points in projective space,
and all the invariants of the theory are invariant under the action of the unitary
group. However, unitary action does not preserve sparse coefficient structure.

In this paper, we give the one-distribution of the roots of random sparse
polynomial systems. We also bound the probability that the condition
number of a random sparse polynomial system is large. 

The roots of a sparse polynomial system are known to belong to a certain
toric variety. However, in order to obtain the theorems below, we needed
to endow the toric variety with a certain geometrical structure, as
explained below. The main insight comes from mechanics, and from symplectic
and K\"ahler geometry.

\subsection{Expected Number of Roots}

Let $A$ be an $M \times n$ matrix, with non-negative
integer entries. To the matrix $A$ we associate the
convex polytope $\cvx (A)$ given by the 
convex hull of all the rows, $\{A^\alpha\}_{\alpha\in\{1,\ldots,M\}}$, of $A$:
\[
\cvx(A) \defeq \left\{ \ \sum_{\alpha=1}^M t_{\alpha} A^{\alpha} 
\ :\ 0 \le t_{\alpha} \le 1, \ \sum_{\alpha=1}^M t_{\alpha} = 1 \ \right \}
\subset (\mathbb R^n)^{\vee}
\punct{.}
\]
\par
  Here, we use the notation $X^{\vee}$ to denote the dual of
a vector space $X$. There are deep reasons to write $\cvx(A)$
as a polytope in dual space, as the reader will see later on.
\par
  Assume that $\dim (\cvx(A)) = n$. Then we can associate to
the matrix $A$ the space $\mathcal F_A$ of polynomials with
support contained in $\{ A^{\alpha} : 1 \le \alpha \le M\}$. This is a 
linear space, and there are many reasonable choices of an
inner product in $\mathcal F_A$. 
\par
  Let $C$ be a diagonal positive definite $M \times M$ matrix. Its inverse
$C^{-1}$ is also a diagonal positive definite $M \times M$ matrix. This
inverse matrix defines the inner product:
\[
  \langle z^{A^{\alpha}}, z^{A^{\beta}} \rangle_{\scriptscriptstyle {C^{-1}}} = (C^{-1})_{\alpha, \beta}
\punct{.}
\]
\par
  The matrix $C$ will be called the {\em variance matrix}.
This terminology arises when we consider random normal 
polynomials in $\mathcal F_A$ with variance $C_{\alpha \alpha}$ for the $\alpha$--th
coefficient.
We will refer to these randomly generated functions as {\em random normal 
polynomials}, for short.  
\medskip
\par
We will also produce several objects associated to the 
matrix $A$ (and to the variance matrix $C$). The most 
important one for this paper will be a K\"ahler manifold 
$(\toric{n}, \omega_A, J)$. This 
manifold is a natural ``phase space'' for the roots
of polynomial systems with support in $A$. It is {\em the} natural
phase space for the roots of systems of random normal polynomials
in $(\mathcal F_A, \langle\cdot,\cdot\rangle_{\scriptscriptstyle {C^{-1}}})$.
\par
More explicitly, let $\toric{n} \defeq 
\mathbb C^n \pmod {2 \pi \sqrt{-1} \,\mathbb Z^n}$ (which, as a real 
manifold, happens to be an $n$-fold product of cylinders). Let 
$\exp: \toric{n} \rightarrow (\mathbb C^*)^n$ denote
coordinatewise exponentiation. Then we will look at the preimages
of the roots of a polynomial system by $\exp$. 
We leave out roots
that have one coordinate equal to zero and roots at infinity.
The differential 2--form $\omega_A$ corresponds to the pull-back
of the canonical 2--form in a suitable Veronese variety (see
Section~\ref{symplectic}). 
\medskip
\par
Systems where all the polynomials have the same support are
called {\em unmixed}. The general situation ({\em mixed} polynomial
systems), where the polynomials may have different supports,
is of greater practical interest. It is also a much more challenging
situation. We shall consider systems of $n$ polynomials in $n$ variables,
each polynomial in some inner product 
space of the form $(\mathcal F_{A_i}, \langle\cdot,\cdot\rangle_{\scriptscriptstyle {C_i^{-1}}})$
(where $i=1, \cdots, n$ and each $A_i$ and each $C_i$ are as
above). 
\par
In this realm, a mathematical object 
(that we may call a {\em mixed manifold})
seems to arise naturally. A mixed manifold is an $(n+2)$--tuple
$(\toric{n}, \omega_{A_1}, \cdots, \omega_{A_n}, J)$ where for each $i$,
$(\toric{n}, \omega_{A_i}, J)$ is a K\"ahler manifold. Mixed manifolds
do {\bf not} have a natural canonical Hermitian structure. They
have $n$ equally important Hermitian structures. However, they have one
natural volume element, the {\em mixed volume form}, given by
\[
d{\toric{n}} = \frac{(-1)^{n(n-1)/2}}{n!} \ 
\omega_{A_1} \wedge \cdots \wedge \omega_{A_n}
\punct{.}
\]
\par
As explained in~\cite{GROMOV},
the volume of $\toric{n}$ relative to the mixed volume form 
is (up to a constant) the mixed volume of the $n$--tuple of polytopes
$(\cvx(A_1), \cdots, \cvx(A_n))$.
\medskip
\par
We extend the famous result by Bernshtein~\cite{BERNSHTEIN} 
on the number of roots of mixed systems of polynomials as follows:

\begin{theorem}\label{GEN-BERNSHTEIN}
 Let $A_1, \cdots, A_n$ and $C_1, \cdots, C_n$ be as above.
For each $i=1,\cdots,n$, let $f_i$ be an (independently 
distributed) normal random polynomial in $(\mathcal F_{A_i}, 
\langle\cdot,\cdot\rangle_{\scriptscriptstyle {C_i^{-1}}})$. Let $U$ be a measurable region of
$\toric{n}$. Then, the expected number of roots of the 
polynomial system $f(z) = 0$ in $\exp U \subseteq (\mathbb C^*)^n$
is
\[
 \frac{n!}{\pi^n}\ \int_U d\toric{n}
\punct{.}
\] 
\end{theorem}

\begin{ex}
  When each $f_i$ is dense with a variance matrix $C_i$
of the form:
\[
C_i = \diag \left(
\frac{\deg f_i !}{I_1 ! I_2 ! \cdots, I_{n} ! (\deg f_i - \sum_{j=1}^n I_j) !}
\right)
\punct{,}\]
the volume element 
$d\toric{n}$ becomes the B\'ezout number $\prod \deg f_i$ times 
the pull-back to \toric{n} of the Fubini-Study metric. 
We thus recover Shub and Smale's stochastic real version of 
B\'ezout's Theorem \cite{BEZ2}. \qed 
\end{ex} 
\par
The general
unmixed case ($A_1 = \cdots = A_n$, $C_1 = \cdots = C_n$) is
a particular case of Theorem~8.1 in~\cite{EDELMAN-KOSTLAN}.
This is the only overlap, since neither theorem generalizes the
other.
\par
  On the other hand, when one sets $U = \toric{n}$, one
recovers Bernshtein's first theorem. The quantity
$\pi^{-n} \int_{\toric{n}} d\toric{n}$ is precisely 
the {\em mixed volume} of polytopes $A_1, \cdots, A_n$
(see~\cite{SANGWINEYAGER} for the classical definition of
Mixed Volume and main properties).
\par
A version of Theorem~\ref{GEN-BERNSHTEIN} was known to
Kazarnovskii~\cite[p.~351]{KAZARNOVSKII} and Khovanskii. 
In~\cite{KAZARNOVSKII}, the supports $A_i$ are allowed
to have complex exponents.
However, uniform variance ($C_i = I$) is assumed. 
His method may imply this special case
of Theorem~\ref{GEN-BERNSHTEIN}, but the indications
given in~\cite{KAZARNOVSKII} were insufficient for us to
reconstruct a proof.
\par
The idea of working with roots of polynomial systems
in logarithmic coordinates seems to be extremely 
classical, yet it gives rise to interesting and
surprising connections (see the discussions in~\cite{GRAEFFE, 
TANGRA,VIRO}).

\medskip
\par

\subsection{The Condition Number}
\label{intro:cond}

  Let $\mathcal F = \mathcal F_{A_1} \times \cdots \times \mathcal F_{A_n}$,
and let $f \in \mathcal F$. A {\em root} of $f$ will be represented
by some $p + q \sqrt{-1} \in \toric{n}$. (Properly speaking, the
root of $f$ is $\exp (p+q\sqrt{-1})$).
\par
  In this discussion, we assume that the ``root'' $p + q \sqrt{-1}$ is
non-degenerate. This means that the derivative of the
{\em evaluation map}
\[
\begin{array}{lrcl}
\mathit{ev}: & \mathcal F \times \toric{n} & \rightarrow & \mathbb C^n \\
& (f, p+q\sqrt{-1}) & \mapsto & (f \circ \exp)(p+q \sqrt{-1})
\end{array}
\]
with respect to the variable in \toric{n}
at the point $p+q\sqrt{-1}$
has rank $2n$. We are then in the situation of the implicit function
theorem, and there is (locally) a smooth function
$G: \mathcal F \rightarrow \toric{n}$ such 
that for $\hat f$ in a neighborhood of $f$, we have
$\mathit{ev}(\hat f, G(\hat f)) \equiv 0$ and $G(f)=p+q\sqrt{-1}$. 
\par
The condition number of $f$ at $(p+q\sqrt{-1})$ is usually defined as
\[
\mymu(f; p+q\sqrt{-1}) = \| DG_{f} \|
\punct{.}
\]
\par
This definition is sensitive to the norm used in the space of 
linear maps between tangent spaces $L(T_f \mathcal F, T_{(p,q)}\toric{n})$.
In general, one would like to use an operator norm, related to some natural
Hermitian or Riemannian structure on $\mathcal F$ and $\toric{n}$.
\par
In the previous Section, we already defined an inner product in each 
coordinate subspace $\mathcal F_{A_i}$, given by the variance matrix
$C_i$. Since the evaluation function is homogeneous in each coordinate,
it makes sense to projectivize each of the coordinate spaces 
$\mathcal F_{A_i}$ (with respect to the inner product $\langle\cdot,
\cdot\rangle_{\scriptscriptstyle {C_i^{-1}}}$).
Alternatively, we can use the Fubini-Study metric in each of the $\mathcal F_{A_i}$'s.
By doing so, we are endowing $\mathcal F$ with a Fubini-like 
metric that is scaling-invariant. We will treat $\mathcal F$ as a multiprojective
space, and write $\mathbb P(\mathcal F)$ for $\mathbb P (\mathcal F_{A_1}) \times
\cdots \times \mathbb P (\mathcal F_{A_n})$.
\par
Another useful metric in $\mathbb P (\mathcal F)$ is given by
\[
d_{\mathbb P} (f,g)^2 \defeq \sum_{i=1}^n \left( \min_{\lambda \in \mathbb C^*} 
\frac{ \|f^i - \lambda g^i\| }{\|f^i\|} \right)^2
\punct{.}
\]
\par
Each of the terms in the sum above corresponds to the square of the
sine of the Fubini (or angular) distance between $f^i$ and $g^i$.
Therefore, $d_{\mathbb P}$ is never larger than the Hermitian distance 
between points in $\cF$, but is a correct first-order aproximation of the 
distance when $g \rightarrow f$ in $\mathbb P (\mathcal F)$. (Compare 
with~\cite[Ch.~12]{BCSS}).
\medskip
\par
While $\mathcal F$ admits a natural Hermitian structure, the solution-space
\toric{n} admits $n$ possibly different Hermitian structures, corresponding
to each of the K\"ahler forms $\omega_{A_i}$. 
\par
In order to elucidate what the natural definition of a condition
number for mixed systems of polynomials is, we will interpret the 
condition number as the  
inverse of the distance to the {\em discriminant locus}. Given
$p+q\sqrt{-1} \in \toric{n}$, we set:
\[
\mathcal F_{(p,q)} = \{ f \in \mathcal F: \mathit{ev}(f;(p,q))=0 \}
\]
and we set $\Sigma_{(p,q)}$ as the space of degenerate polynomial systems
in $\mathcal F_{(p,q)}$. Since the fiber $\mathcal F_{(p,q)}$ inherits
the metric structure of $\mathcal F$, we can speak of the 
distance to the discriminant locus along a fiber. In this setting,
Theorem~3 in~\cite[p.~234]{BCSS} becomes:
\begin{theorem}[Condition number theorem]~\label{condnumber}
  Under the notations above, if $(p,q)$ is a non-degenerate root of $f$,
\[
\max_{ \|\dot f\| \le 1}
\min_{i}
\| DG_f \dot f \|_{A_i}
\le
\frac{1}{d_{\mathbb P}(f, \Sigma_{(p,q)})}
\le
\max_{ \|\dot f\| \le 1}
\max_{i}
\| DG_f \dot f \|_{A_i}
\punct{.}
\]
\end{theorem} 
There are two interesting particular cases. First of
all, if $A_1 = \cdots = A_n$ and $C_1 = \cdots = C_n$,
we obtain an equality:

\begin{corollary}[Condition number theorem for unmixed systems]
\ \newline
  Let $A_1 = \cdots = A_n$ and $C_1 = \cdots = C_n$, then
  under the hypotheses of Theorem~\ref{condnumber},
\[
\mymu(f; (p,q))
\defeq
\max_{ \|\dot f\| \le 1}
\min_{i}
\| DG_f \dot f \|_{A_i}
=
\max_{i}
\max_{ \|\dot f\| \le 1}
\| DG_f \dot f \|_{A_i}
=
\frac{1}{d_{\mathbb P}(f, \Sigma_{(p,q)})}
\punct{.}
\]
\end{corollary}

We can also obtain a version of 
Shub and Smale's condition number theorem (Theorem~3 in~\cite[p.~243]{BCSS})
for dense systems as a particular case, once we choose the
correct variance matrices:
\begin{corollary}[Condition number theorem for dense systems]
  \ \newline Let $d_1, \cdots, d_n$ be positive integers, and
let $A_i$ be the $n$-columns matrix having all possible 
rows with non-negative entries adding up to at most $d_i$.
Let
\[
C_i = 
\frac{1}{d_i} \
\diag \left( 
\frac{d_i !}{ 
(A_i)^{\alpha}_1 !
(A_i)^{\alpha}_2 !
\cdots
(A_i)^{\alpha}_n !
(d_i - \sum_{j=1}^n (A_i)^{\alpha}_j) !
} \right)
\punct{.}
\]
Then,
\[
\mymu(f; (p,q))
\defeq
\max_{ \|\dot f\| \le 1}
\min_{i}
\| DG_f \dot f \|_{A_i}
=
\max_{i}
\max_{ \|\dot f\| \le 1}
\| DG_f \dot f \|_{A_i}
=
\frac{1}{d_{\mathbb P}(f, \Sigma_{(p,q)})}
\punct{.}
\]
\end{corollary}

  The factor $\frac{1}{d_i}$ in the definition of
the variance matrix $C_i$ 
corresponds to the factor $\sqrt{d_i}$ in the definition
of the normalized condition number in~\cite[p.~233]{BCSS}. 
It scales the K\"ahler forms $\omega_{A_i}$ so that they
are equal (see Remark~\ref{recoverkostlan} p.~\pageref{recoverkostlan}
below).
\medskip
\par
In the general mixed case, we would like to interpret the 
two ``minmax'' bounds as condition numbers related to
some natural Hermitian or Finslerian structures on
\toric{n}. See Section~\ref{finsler} for a discussion.
\par
Theorem~\ref{condnumber} is very similar to Theorem~D in
\cite{DEDIEU-PARK-CITY}, but the philosophy here is
radically different. Instead of changing the metric in
the fiber $\mathcal F_{(p,q)}$, we consider the inner
product in $\mathcal F$ as the starting point of our 
investigation. Theorem~\ref{condnumber} gives us some
insight about reasonable metric structures
in $\mathcal T^n$.
\medskip
\par

As in Theorem~\ref{GEN-BERNSHTEIN}, let $U$ be
a measurable set of \toric{n}. In view of Theorem~\ref{condnumber},
we define a restricted condition number (with respect to
$U$) by:
\[
\mymu (f; U)
\defeq
\frac{1}{\min_{(p,q) \in U} d_{\mathbb P}(f, \Sigma_{(p,q)})}
\]
where the distance $d_{\mathbb P}$ is taken along the fiber
$\mathcal F_{(p,q)} = \{ f : (f\circ \exp) (p+q \sqrt{-1}) = 0 \}$.

Although we do not know in general how to bound the expected value of
$\mymu(f; \toric{n})$, we can give a convenient bound for
$\mymu(f; U)$ whenever $U$ is compact and in some cases where
$U$ is not compact. 

The group $GL(n)$ acts on $T_{(p,q)}\mathcal T^n$ by sending $(\dot p,\dot q)$
into $(L\dot p,L\dot q)$, for any $L \in GL(n)$. 
In more intrinsic terms, $J$ and the $GL(n)$-action
commute. With this convention, we can define an intrinsic invariant of
the mixed structure $(\mathcal T^n, \omega_{A_1}, \cdots, \omega_{A_n}, J)$:

\begin{definition}\label{intrinsic}
  The {\em mixed dilation} of the tuple $(\omega_{A_1}, \cdots, \omega_{A_n})$
is:
\[
\kappa(\omega_{A_1},\cdots,\omega_{A_n}; (p,q)) \defeq
\min_{L \in GL(n)}
\max_i
\frac{
  \max_{\|u\|=1} (\omega_{A_i})_{(p,q)} (Lu,JLu)
}{
  \min_{\|u\|=1} (\omega_{A_i})_{(p,q)} (Lu,JLu)
}
\punct{.}
\]
\end{definition}

Given a set $U$, we define:
\[
\kappa_{U} \defeq \sup_{(p,q) \in U} \kappa(\omega_{A_1},\cdots,\omega_{A_n}; (p,q))
\punct{,}
\]
provided the supremum exists, and $\kappa_U = \infty$ otherwise.

We will bound the expected number of roots with condition number
$\mymu > \eps^{-1}$ on $U$ in terms of the mixed volume form,
the mixed dilation $\kappa_U$ and the expected 
number of ill-conditioned roots
in the {\em linear case}. The linear case corresponds to polytopes
and variances below:

\begin{align*}
A^{\linear}_i &=
\left[
\begin{matrix}
0 & \cdots & 0 \\
1 & & \\
  & \ddots & \\
  & & 1  
\end{matrix}
\right]
&
C^{\linear}_i =
\left[
\begin{matrix}
1 \\
  & 1 \\
  &   & \ddots & \\
  &   &        & 1  
\end{matrix}
\right]
\end{align*}

\begin{theorem}[Expected value of the condition number]
\label{mixed-cond} 
Let $\nu^{\linear}(n,\eps)$ be the probability that a random 
$n$--variate linear complex polynomial has condition number
larger than $\eps^{-1}$. Let $\nu^{A}(U,\eps)$ be the
probability that $\mymu(f,U) > \eps^{-1}$ for a 
normal random polynomial system
$f$ with supports $A_1,\cdots,A_n$ and variance $C_1,\cdots,C_n$.

Then,
\[
\nu^A(U,\eps)
\le
\frac{ \int_U \bigwedge \omega_{A_i} }
{ \int_U \bigwedge \omega_{A^{\linear}_i}}
\
\nu^{\linear} (n, \sqrt{\kappa_U} \eps)
\punct{.}
\]
\end{theorem}

There are a few situations where we can assert that $\kappa_U = 1$.
For instance,

\begin{corollary}\label{unmixed:complex}
  Under the hypotheses of Theorem~\ref{mixed-cond}, if
$A = A_1 = \cdots = A_n$ and $C=C_1 = \cdots = C_n$, 
then
\[
\nu^A(U,\eps)
\le
\vol(U)
\ \nu^{\linear} (n, \eps)
\punct{.}
\]
\end{corollary}

The dense case (Theorem~1 p.~237 in~\cite{BCSS}) is also a consequence
of Theorem~\ref{mixed-cond}.

\begin{remark}
   We interpret $\nu^{\linear}(n,\eps)$ as the probability that
a random linear polynomial $f$ is at multiprojective 
distance less than $\eps$ from the discriminant variety
$\Sigma_{(p,q)}$.
   Let $g \in \Sigma_{(p,q)}$ be such that the following minimum
is attained:
\[
d_{\mathbb P} (f, \Sigma_{(p,q)}) ^2 = 
\inf_{\substack{g \in \Sigma_{(p,q)}\\ \lambda \in (\mathbb C^*)^n}}
\sum_{i=1}^n \frac{\|f^i - \lambda_i g^i\|^2}{\|f^i\|^2}
\punct{.}
\]
Without loss of generality, we may scale $g$ such that $\lambda_1 = 
\cdots = \lambda_n = 0$. In that case, 
\[
d_{\mathbb P} (f, \Sigma_{(p,q)}) ^2 = 
\sum_{i=1}^n \frac{\|f^i - g^i\|^2}{\|f^i\|^2}
\ge
\frac{ \sum_{i=1}^n \|f^i - g^i\|^2}{\sum_{i=1}^n \|f^i\|^2}
\punct{.}
\]
The right hand term is the projective distance to the
discriminant variety along the fiber, in the sense of~\cite{BCSS}.
Since we are in the linear case, this may be interpreted as
the inverse of the condition number of $f$ in the sense
of~\cite[Prop.~4 and Remark~2 p.~250]{BCSS}. 
\par
Recall that each $f^i$ is an independent random normal linear polynomial
of degree 1, and that $C_i$ is the identity. Therefore, each
$f^i_{\alpha}$ is an i.i.d. Gaussian variable. If we look at the system
$f$ as a random variable in $\mathbb P^{n(n+1)-1}$, then we obtain 
the same probability distribution as in~\cite{BCSS}.
Then, using Theorem~6 p.~254 {\em ibid}, we deduce that
\[
\nu^{\linear}(n,\eps) \le \frac{n^3 (n+1) \Gamma(n^2+n)}{\Gamma(n^2+n-2)} \eps^4
\punct{. \qed}
\]
\end{remark}

\subsection{Real Polynomials}

Shub and Smale showed in~\cite{BEZ2}
that the expected number of real roots, in the dense
case (with unitarily invariant probability measure)
is exactly the square root of the expected number of roots.
\par
Unfortunately, this result seems to be very
hard to generalize to the unmixed case. Under
certain conditions, explicit formul\ae\ for the unmixed case
are available~\cite{ROJAS-PARK-CITY}. Also,
less explicit bounds for the multi-homogeneous case were
given by~\cite{MCLENNAN}.
\par
Here, we will give a very coarse estimate in
terms of the square root of the mixed volume:

\begin{theorem} \label{expected-real}
  Let $U$ be a measurable set in $\mathbb R^n$,
with total Lebesgue volume $\lambda(U)$. Let $A_1,\cdots,A_n$
and $C_1,\cdots,C_n$ be as above. Let
$f$ be a normal random real polynomial system.
Then the average number of real roots of $f$
in $\exp U \subset \mathbb (R^+_*)^n$ is 
bounded above by
\[
(4\pi^2)^{-n/2}
\sqrt{\lambda(U)}
\sqrt{
\int_{\substack{(p,q)\in \toric{n}\\ p\in U}} n! d\toric{n}
}
\punct{.}
\]
\end{theorem}

This is of interest when $n$ and $U$ are fixed.
In that case, the expected number of positive
real roots (hence of real roots) grows as
the square root of the mixed volume.

It is somewhat easier to investigate real random
polynomials in the unmixed case. 

Let $\nu_{\mathbb R}(n,\eps)$
be the probability that a
linear random real polynomial 
has condition number
larger than $\eps^{-1}$.  

\begin{theorem}\label{unmixed:real}
  Let $A = A_1 = \cdots = A_n$ and $C=C_1 = \cdots = C_n$. 
Let $U \subseteq \mathbb R^n$ be measurable. Let $f$
be a normal random real polynomial system. Then,
\[
\mathrm{Prob} \left[ \mymu(f,U) > \eps^{-1} \right]
\le
E(U) \ \nu_{\mathbb R}(n,\eps)
\]
where $E(U)$ is the expected number of real roots on $U$.
\end{theorem}

Notice that $E(U)$ depends on $C$. Even if we make
$U = \mathbb R^n$, we may still obtain a bound depending
on $C$.

\subsection{Acknowledgements}

Steve Smale provided valuable inspiration for us to develop 
a theory similar to~\cite{BEZ1,BEZ2,BEZ3,BEZ4,BEZ5}
for sparse polynomial systems. He also provided
examples on how to eliminate the dependency upon unitary
invariance in the dense case.
\par
  The paper by Gromov~\cite{GROMOV} was of foremost importance to
this research. To the best of our knowledge, ~\cite{GROMOV} is the only clear
exposition available of mixed volume in terms of a wedge of differential
forms. We thank Mike Shub for pointing out that reference, and for
many suggestions.
\par
  We would like to thank Jean-Pierre Dedieu for sharing his thoughts
with us on Newton iteration in Riemannian and quotient manifolds.  
\par
  Also, we would like to thank Felipe Cucker, Alicia Dickenstein,
Ioannis Emiris, Askold Khovanskii, Eric Kostlan, T.Y. Li, 
Martin Sombra and Jorge P.\ Zubelli for their suggestions and support.
\medskip
\par
  This paper was written while G.M. was visiting the Liu Bie Ju Center
for Mathematics at the City University of Hong Kong. He wishes to thank
CityU for the generous support.

\section{Symplectic Geometry and Polynomial Systems}
\label{symplectic}

\subsection{About Symplectic Geometry}

\begin{definition}[Symplectic structure]
  Let $M$ be a manifold. A $2$--form on $M$ is said to be
non-degenerate if and only if for all $x \in M$, the only
vector $u \in T_x M$ such that for all $v \in T_x M$,
$\omega_x (u,v) = 0$
is the zero vector.
\par
  A {\em symplectic form} on $M$ is a closed, non-degenerate $2$--form
$\omega$ on $M$. In that case, $(M, \omega)$ is said to be a symplectic
manifold.
\end{definition}

\begin{definition}[Complex structure]
  Let $M$ be a complex manifold. (We assume that $M$ is given with a
certain maximal holomorphic atlas). If $X: U \subset \mathbb C^n \rightarrow M$ 
is a chart of $M$, and $p = X(z) \in M$,
then we define $J_p: T_pM \rightarrow T_pM$ such that the following
diagram commutes:
\[
\begin{CD}
 T_pM @>{J_p}>> T_pM \\
 @A{DX_z}AA  @A{DX_z}AA \\
 T_z \mathbb C^n @>{\text{Mult. by $\sqrt{-1}$}}>> T_z \mathbb C^n \\
\end{CD}
\punct{.}
\]
\par
 This is well-defined for each $p$ in $M$. Indeed, if two charts $X$ and $Y$ overlap
at $p$, then $Y \circ X^{-1}: \mathbb C^n \rightarrow \mathbb C^n$ is holomorphic so 
its derivative exists and commutes with multiplication by $\sqrt{-1}$.
\par
The map
\[
\begin{array}{lrcl}
 J: & TM & \rightarrow & TM \\
 & (p,\dot p) & \mapsto & (p, J_p \dot p)
\end{array} 
\]
is called the {\em canonical complex structure} of $M$. (The complex structure may
depend on the holomorphic atlas. We assume a canonical holomorphic
atlas of $M$ is given). Note that $-J^2$ is the identity on $TM$.
\end{definition}

\begin{definition}[K\"ahler manifolds]
  Let $M$ be a complex manifold, with complex structure $J$. Let $\omega$ be
a symplectic form on $M$ (considered as a real manifold). Then $\omega$ is
said to be a {\em $(1,1)$--form} if and only if $J^* \omega = \omega$. A 
$(1,1)$ form $\omega$ corresponds to a symmetric form $u,v \mapsto \omega(u,Jv)$.
We say that $\omega$ is strictly positive if and only if the corresponding
symmetric form is positive definite for all $p \in M$. Therefore, a
strictly positive $(1,1)$--form defines a Riemann structure on $M$. Also,
we can recover an Hermitian structure on $M$ by setting 
$\langle u,v \rangle = \omega(u,Jv) + \sqrt{-1} \omega(u,v)$. 
\par
  The triple $(M, \omega, J)$ is said to be a K\"ahler manifold when
$M$, $\omega$ and $J$ are as above.
\end{definition}

\begin{example}[$\mathbb C^M$] 
  We identify $\mathbb C^{M}$ to $\mathbb R^{2M}$ and use coordinates
$Z^i = X^i + \sqrt{-1} Y^i$. The {\em canonical} $2$--form 
$\omega_Z = \sum_{i=1}^M dX_i \wedge dY_i$ makes $\mathbb C^{M}$ into a
symplectic manifold.
\par
  The natural complex structure $J$ is just the multiplication by ${\sqrt{-1}}$.
The triple $(\mathbb C^M, \omega_Z, J)$ is a K\"ahler manifold. \qed 
\end{example}

\begin{example}[Projective space]\label{ex:proj}
  Projective space $\mathbb P^{M-1}$ admits a {\em canonical} $2$--form
defined as follows. Let $Z = (Z^1, \cdots, Z^M) \in (\mathbb C^{M})^*$,
and let $[Z] = (Z^1 : \cdots : Z^M) \in \mathbb P^{M-1}$ be the corresponding 
point in $\mathbb P^{M-1}$.
  The tangent space $T_{[Z]} \mathbb P^{M-1}$ may be modelled by $Z^{\perp}
\subset T_Z \mathbb C^M$. 
  Then we can define a two-form on $\mathbb P^{M-1}$ by setting:
\[
\omega_{[Z]} (u,v) = \|Z\|^{-2} \omega_Z (u,v) \punct{,}
\]
  where it is assumed that $u$ and $v$ are orthogonal to $Z$. The 
latter assumption tends to be quite inconvenient, and most people 
prefer to pull $\omega_{[Z]}$ back to $\mathbb C^M$ by the canonical 
projection $\pi: Z \mapsto [Z]$.
  It is standard to write the pull-back $\tau = \pi^* \omega_{[Z]}$
as:
\[
\tau_Z = - \frac{1}{2}d J^* d \ \frac{1}{2}\log \|Z\|^2 \punct{,}
\]
using the notation $d\eta = \sum_i \frac{\partial \eta}{p_i} \wedge dp_i
+ \frac{\partial \eta}{q_i} \wedge dq_i$, and where $J^*$ denotes the
pull-back by $J$.
\par
  Projective space also inherits the complex structure from $\mathbb C^M$.
Then $\omega_{[Z]}$ is a strictly positive $(1,1)$--form. The corresponding
metric is called {\em Fubini-Study} metric in $\mathbb C^M$ or 
$\mathbb C^{M-1}$. \qed 
\end{example}

\begin{remark}
  Some authors prefer to write 
$\sqrt{-1} \partial \bar \partial$
instead of $- \frac{1}{2} d J^* d$. The following notation
is assumed:
$\partial \eta = \sum_i \frac{\partial \eta}{Z_i} \wedge dZ_i$
and $\bar \partial \eta = \sum_i \frac{\partial \eta}{\bar Z_i} \wedge d\bar Z_i$.
Then they write $\tau_Z$ as:
\[
\tau_Z = \frac{\sqrt{-1}}{2}
\left(
\frac{ \sum_i dZ_i \wedge d\bar Z_i }{\|Z\|^2}
-
\frac{ \sum_i Z_i d\bar Z_i \wedge \sum_j \bar Z_j dZ_j}{\|Z\|^4}
\right)
\punct{. \qed}
\]
\end{remark}

\begin{example} Let $A$ be an $M \times n$ matrix with non-negative
integer entries, and we associate every row $A^{\alpha}$ of $A$ to
the monomial $z^{A^{\alpha}} = z_1 ^{A^{\alpha}_1} \cdots z_n ^{A^{\alpha}_n}$. 
We also assume (as in the introduction)
that the corresponding polytope (the convex hull of all the rows)
is $n$-dimensional. Also, as in the introduction, let $C$ be a diagonal
positive-definite matrix (that we called the variance matrix). The
variance matrix was the matrix of the inner product in $\mathcal F_A$.
Let $C^{1/2}$ be
the diagonal positive-definite matrix such that $C=C^{1/2} C^{1/2}$.
The right-multiplication of
some $f \in \mathcal F_A$ by $C^{-1/2}$ makes the inner product canonical.
The left-multiplication by $C^{1/2}$ is the pull-back of this operation
in dual-space $\mathcal F_A^{\vee}$.
\medskip
\par
We define the map $\hat V_A$ from $\mathbb C^n$ into $\mathbb C^{M}$:
\[
\hat V_A: z \mapsto 
C^{1/2}
\left[
\begin{matrix} z^{A^1} \\ \vdots \\ z^{A^M} \end{matrix} \right]
\punct{.}
\]
\par
Because $C^{1/2}$ is diagonal, 
$\| \hat V_A(z) \|$ is invariant by the natural action 
$z_i \mapsto z_i e^{\theta_i \sqrt{-1}}$ of
$(S^1)^n$ on the variable $z$. Moreover, we still have
the pairing $f(z) = (f \cdot C^{-1/2}) \cdot \hat V_A(z)$.
The variable $(f \cdot C^{-1/2})$ is corresponds to
$M_i$ i.i.d. Gaussian variables.
\medskip
\par
We can also compose with the projection into projective space,
$V_A = \pi \circ \hat V_A: \mathbb C^n \rightarrow \mathbb P^{M-1}$,
\[
V_A: z \mapsto [\hat V_A(z)]
\punct{.}
\] 
\par
When $C$ is the identity, the Zariski closure of the
image of $V_A$ is called the {\em Veronese variety}. The map $V_A$ is
called the {\em Veronese embedding}. Notice that $V_A$ is not defined
for certain values of $z$, like $z=0$. Those values are called the
{\em exceptional set}. The exceptional set is contained in the union
of the planes $z_i=0$. 
\par
There is a natural symplectic structure on the closure of the image of
$V_A$, given by the restriction of the Fubini-Study 2--form. We will
see below (Lemma~\ref{dimnondeg}) 
that $DV_A$ has rank $n$ for $z \in (\mathbb C^*)^n$, because
the polytope of $A$ has dimension $n$. 
Thus, we can pull-back this structure to $\mathbb (C^*)^n$ by:
\[
\Omega_A = V_A ^* \tau
\punct{.}
\]
\par
Also, we can pull back the complex structure of $\mathbb P^{M-1}$,
so that $\Omega_A$ becomes a strictly positive $(1,1)$-form. 
\par
Therefore, the matrix $A$ defines a K\"ahler manifold
$(\mathbb (C^*)^n, \Omega_A, J)$. \qed 
\end{example}

\begin{example} Let $\toric{n} = \mathbb C^n \pmod {2 \pi \sqrt{-1} \,\mathbb Z^n}$.
We will use coordinates
$p + q \sqrt{-1}$ for $\toric{n}$, where $p \in \mathbb R^n$ and
$q \in \mathbb T^n = \mathbb R^n \pmod {2 \pi \,\mathbb Z^n}$.
\par
Given $(p,q) \in \toric{n}$, we define $\exp (p,q) = \left(\cdots,e^{p+q\sqrt{-1}},\cdots\right) 
\in (\mathbb C^*)^n$. 
\par
For any matrix $A$ as in the previous example, we can pull-back the K\"ahler
structure of $((\mathbb C^*)^n, \Omega_A, J)$ to obtain another K\"ahler
manifold $(\toric{n}, \omega_A, J)$.
(Actually, it is the same object in
logarithmic coordinates, minus points at ``infinity''.) 
An equivalent definition is to pull back the K\"ahler structure
of the Veronese variety by
\[
\hat v_A \defeq \hat V_A \circ \exp 
\punct{. \qed}
\]
\end{example}

\begin{remark}
  The Fubini-Study metric on $\mathbb C^M$ was constructed by applying the
operator $-\frac{1}{2}dJ^*d$ to a certain convex function
(in our case, $\frac{1}{2}\log \|Z\|^2$). This is a general standard way to construct
K\"ahler structures. In~\cite{GROMOV}, it is explained how to associate a
(non-unique) convex function to any convex body, thus producing an associated
K\"ahler metric. \qed 
\end{remark}

\begin{remark} Now a little bit of magic... 
$\omega_A = \hat v ^* \tau = 
\hat v ^* (- \frac{1}{2}dJ^*d) g$, where $g: Z \mapsto \frac{1}{2} \log \|Z\|^2$.
Both $d$ and $J$ commute with pull-back, so
\begin{eqnarray*}
\omega_A &=& 
-  
\hat v ^* 
(\frac{1}{2}dJ^*d) g \\
&=&
-  
(\frac{1}{2}d J^* d) \hat v_A^* g
\\
&=&
-  
(\frac{1}{2} d J^* d) (g \circ \hat v_A) \ \text{\qed} 
\end{eqnarray*}
\end{remark}

\begin{remark} \label{recoverkostlan} 
The same is true for $((\mathbb C^*)^n, \Omega_A, J)$.
A particular case should be mentioned here. 
Unitary invariance played an important role in~\cite{BEZ1,BEZ2,BEZ3,BEZ4,BEZ5} and
in~\cite{BCSS}. Let us now recover that invariance for dense 
polynomials.
\par
Suppose that the rows of our matrix $A$
are the exponent vectors of all possible monomials of degree 
exactly $D$ in $n+1$ variables.
Let $A^{\alpha} = [I_1, \cdots, I_{n+1}]$. We set 
$C = \diag\left(\frac{D!}{I_1!I_2! \cdots, I_{n+1}}\right)$. Then,
\[
g \circ \hat V_A = \frac{1}{2}\log \| \hat V_A \|^2 
= \frac{1}{2} \log \| z_1, \cdots, z_{n+1} \|^{2D}
= D \frac{1}{2}\log \|z_1, \cdots, z_{n+1}\|^2
\punct{.}
\]
So $\Omega_A$ is a multiple of the Fubini-Study metric, and 
we can actually extend $\Omega_A$ to $\mathbb C^{n+1}_*$.
\par
 Let $\tilde f = f C^{-1/2} \in (\mathbb C^M)^{\vee}$ represent the 
polynomial $z \mapsto f(z) = \tilde f \hat V_A(z)$. Then
evaluation corresponds to the pairing $(\tilde f,z) \mapsto \tilde f \cdot V_A(z)$. 
The action of $U(n)$ on $\mathbb C^{n+1}$ induces 
an action on $(\mathbb C^M)^{\vee}$. 
All these actions are unitary, and the Hermitian structure
of the space of polynomials (in the coordinates above) is
invariant under such actions. \qed 
\end{remark}

For the record, we state explicit formul\ae\ for several of the
invariants associated to the K\"ahler manifold
$(\toric{n}, \omega_A, J)$. First of all, the function
$g_A = g \,\circ \,\hat v_A$ is precisely:

\formula{\label{211}The canonical Integral $g_A$ (or {\em K\"ahler potential})
of the convex set associated to $A$}
{ g_A(p) = \frac{1}{2} \log 
\left(
\left( \exp (A \cdot p) \right) ^T
C
\left( \exp (A \cdot p) \right) \right) }

The terminology {\em integral} is borrowed from mechanics,
and its refers to the invariance of $g_A$ by $\mathbb T^n$-action
(see appendix ~\ref{momentum} for more analogies). Also, the
gradient of $g_A$ is called the {\em momentum map}. Recall that 
the Veronese
embedding takes values in projective space. We will use the
following notation:
$v_A(p) = \hat v_A(p) / \| \hat v_A(p) \|$. 
This is independent
of the representative of equivalence class $v_A(p)$. 
Now, let $v_A(p)^2$ mean
coordinatewise squaring and $v_A(p)^{2T}$ be the transpose 
of $v_A(p)^2$. The gradient of $g_A$ is then:
\formula{The Momentum Map associated to $A$}
{
\nabla g_A = v_A(p)^{2T} A
}

Since $p \mapsto v_A(p)$ is a well-defined real function,
we may write its derivative as 
\[
Dv_A(p) = P_{v_A(p)}
\diag(v_A(p)) A
\]
where $P_v$ is the projection operator $I-\frac{v v^H}{\|v\|^2}$.

Then the second derivative of $g_A$ is

\formula{\label{secder} Second derivative of $g_A$}
{
D^2 g_A = 2 Dv_A(p)^T Dv_A(p)
}

Using the relation $-\frac{1}{2}dJ^*d g_A = \frac{1}{2}\sum (D^2 g_A)_{ij} dp_i \wedge dq_j$,
one obtains an expression for $\omega_A$:

\formula{\label{omega} The symplectic 2--form associated to $A$:}
{
(\omega_A)_{(p,q)} = \frac{1}{2}\sum_{ij} (D^2 g_A)_{ij} dp_i \wedge dq_j
}

We still have to show that $\omega_A$ is a symplectic form. Clearly,
$\frac{1}{2}d (d J^* d g_A) = \frac{1}{2}d^2 (J^* d g_A) = 0$. The remaining condition to
check is non-degeneracy. In view of formul\ae ~\ref{secder} and
~\ref{omega}, this is a consequence of the following fact:

\begin{lemma} \label{dimnondeg}
Let $A$ be a matrix with non-negative integer entries,
such that $\cvx(A)$ has dimension $n$. Then $(Dv_A)_p$ is injective,
for all $p \in \mathbb R^n$.
\end{lemma}

\begin{proof}
  The conclusion of this Lemma can fail only if there are $p \in \mathbb R^n$
and $u\ne 0$
with $(Dv_A)_p u = 0$. This means that
\[
P_{v_A(p)} \mathrm{diag} (v_A)_p A u = 0
\punct{.}
\]

This can only happen if $\mathrm{diag} (v_A)_p A u$ is in the
space spanned by $(v_A)_p$, or, equivalently, $Au$ is in the
space spanned by $(1,1, \cdots, 1)^T$. This means that all the
rows $a$ of $A$ satisfy $a u = \lambda$ for some $\lambda$. 
Interpreting a row of $A$ as
a vertex of $\cvx A$, this means that $\cvx A$ is contained in
the affine plane $\{ a : a u = \lambda \}$.
\end{proof}

We can also write down the Hermitian structure of \toric{n} as:

\formula{Hermitian structure of \toric{n} associated to $A$:}
{
(\langle u, w \rangle_{\scriptscriptstyle A})_{(p,q)} 
= 
u^H (\frac{1}{2}D^2 g_A)_p w
}

In general, the function $v_A$ goes from \toric{n}
into projective space. Therefore, its derivative
is a mapping

\[
(Dv_A)_{(p,q)}:
T_{(p,q)}\toric{n}
\rightarrow T_{\scriptscriptstyle v_A(p+q \sqrt{-1})} \mathbb P^{M-1}
\simeq \hat v_A(p+q\sqrt{-1}) ^{\perp} \subset \mathbb C^M
\punct{.}
\]

For convenience, we will write this derivative as a mapping into
$\mathbb C^M$, with range $\hat v_A(p+q\sqrt{-1}) ^{\perp}$.
Let $P_v$ be the projection operator
\[
P_v = I - \frac{1}{\|v\|^2} v v^H
\punct{.}
\]
Then,

\formula{\label{deriv} Derivative of $v_A$}
{
(Dv_A)_{\scriptscriptstyle (p,q)} = 
P_{ \hat v_A(p+q \sqrt{-1}) }
\diag \left( \frac{ \hat v_A(p+q \sqrt{-1}) }{\|\hat v_A(p+q \sqrt{-1}\|}
\right) A
}

An immediate consequence of Formula~\ref{deriv} is:

\begin{lemma}\label{lem-orthog}
  Let $f \in \mathcal F_A$ and $(p,q) \in \toric{n}$
be such that $f \cdot \hat v_A(p+q\sqrt{-1})=0$. Then, 
$f \cdot (Dv_A)_{(p,q)} = \frac{1}{\| \hat v_A(p,q)\|} (D\hat v_A)_{(p,q)}$
\end{lemma}

  In other words, when $(f \circ \exp) (p+q\sqrt{-1})$ vanishes,
$Dv_A$ and $D\hat v_A$ are the same up to scaling.

Notice that the Hermitian metric is also 
\[
(\langle u,w\rangle_{\scriptscriptstyle A})_{(p,q)} = u^h Dv_A(p,q)^H Dv_A(p,q) w
\punct{.}
\]

Finally, the volume element associated to $A$ is
\formula{Volume element of $(\toric{n}, \omega_A, J)$}
{
d\toric{n}_A = \det \left( \frac{1}{2} \ D^2 g_A(p)\right) \  
dp_1 \wedge \cdots \wedge dp_n \wedge dq_1 \wedge \cdots \wedge dq_n 
}

\subsection{Toric Actions and the Momentum Map}

The {\em momentum map}, also called {\em moment map}, was 
introduced in its modern formulation by Smale~\cite{SMALE70}
and Souriau~\cite{SOURIAU}. The reader may consult one of
the many textbooks in the subject (such as Abraham and
Marsden~\cite{ABRAHAM-MARSDEN} or McDuff and Salamon
~\cite{MCDUFF-SALAMON}) for a general exposition. 
\par
In appendix~\ref{momentum}, we will explicitly show
what Lie group action $\nabla g_A$ is
the momentum of, and what
the associated Hamiltonian dynamical system is. 
\par
In this Section we instead follow the point of view
of Gromov~\cite{GROMOV}. The main results in this
Section are that

\begin{proposition} \label{convexity}
  The momentum map $\nabla g_A$ maps 
$\toric{n}$ onto the interior of $\cvx(A)$.
When $\nabla g_A$ is restricted to the
real $n$-plane $[q=0] \subset \toric{n}$,
this mapping is a bijection. \qed 
\end{proposition}

  This seems to be a particular case of the 
Atiyah-Guillemin-Sternberg theorem. However, technical
difficulties prevent us from directly applying this result here
(see appendix~\ref{momentum}).

\begin{proposition} \label{volpres}
  The momentum map $\nabla g_A$ is a volume-preserving
map from the manifold $(\toric{n}, \omega_A, J)$
into $\cvx(A)$, up to a constant, in the following
sense: if $U$ is a measurable region of $\cvx(A)$,
then
\[
\vol \left( (\nabla g_A)^{-1}(U) \right) = \pi^{n}\  \vol U
\punct{. }
\]
\end{proposition}

We prove Proposition~\ref{volpres} by first assuming 
Proposition~\ref{convexity}.  

\begin{proof}[Proof of Proposition~\ref{volpres}]
Consider the mapping
\[
\begin{array}{rrcl}
M: & \toric{n} & \rightarrow &  \frac{1}{2}\cvx(A) \times \mathbb T^n \\
& (p,q) & \mapsto & (\frac{1}{2}\nabla g_A(p), q)
\end{array}
\punct{.}
\]

Since we assume $\dim \cvx(A) = n$, we can apply Proposition~\ref{convexity}
and conclude that $M$ is a diffeomorphism. 
\medskip
\par
The pull-back of the canonical symplectic structure in
$\mathbb R^{2n}$ by $M$ is precisely $\omega_A$, because
of Formul\ae~\ref{secder} and~\ref{omega}. 
Diffeomorphisms with that property are called
{\em symplectomorphisms}. Since the volume form of a
symplectic manifold depends only of the canonical 2--form,
symplectomorphisms preserve volume. We compose with a scaling by 
$\frac{1}{2}$ in the first $n$ variables, that divides $\vol U$ by $2^n$, 
and we are done.
\end{proof}
\begin{remark}
  Symplectomorphisms are also known to preserve a few other invariants
such as the symplectic width (see ~\cite{MCDUFF-SALAMON}).
However, symplectomorphisms are not required to 
preserve the complex structure and therefore need not
be isometries. 
\par
  However, it is explained in~\cite{ABREU} how to define a new
complex structure in $\cvx(A) \times \mathbb T^n$ that will make
the map $M$ a K\"ahler isomorphism, hence an isometry. \qed 
\end{remark}

Before proving Proposition~\ref{convexity}, we will need
the following result about convexity. We follow here
Convexity Theorem 1.2 in ~\cite{GROMOV}, attributed to Legendre:
\begin{theorem*}[Legendre]
 If $f$ is convex and of class $\mathcal C^2$ on $\mathbb R^n$, 
then the closure of the image
 $\{ \nabla f _{r} : r \in \mathbb R^n\}$
in $\mathbb (R^n)^{\vee}$ is convex.
\end{theorem*}

\begin{proof}
  Let $L_f$ be the set of covectors $y\in (\mathbb R^n)^{\vee}$
with the property that 
\[
\exists c \in \mathbb R \ \forall x \in \mathbb R^n  f(x) \ge y \cdot x -c
\punct{.}
\]
\par
 Notice that $L_f$ is a convex subset of $(\mathbb R^n)^{\vee}$. Geometrically,
the planes in $L_f$ with $c$ minimal correspond to the {\em envelope}
of the graph of $f$.
\par
 The set $L_f$ contains $\{ \nabla f _{r} : r \in \mathbb R^n\}$:
For any given $r$, we set $c_r = \nabla f_r \cdot r - f(r)$. 
Since $f$ is convex, 
\[
f(x) \ge \nabla f_r \cdot x - c_r  
\punct{.}
\]
\par
 For the converse, assume that there is $y \in L_f$ not in
the closure of $\{ \nabla f_{r} : r \in \mathbb R^n\}$.
Then there is some $\eps > 0$ such that
\[
\forall r \in \mathbb R^n,  \| y - \nabla f_r \| > \eps
\punct{.}
\]
\par
We define the following gradient vector field in $\mathbb R^n$:
\[
\dot x = \frac{ (y - \nabla f_x)^T }{\| y - \nabla f_x\|}
\]
Because the denominator is bounded below by $\eps$, this field
is well-defined and Lispchitz in all of $\mathbb R^n$. Let us fix
an arbitrary initial condition $x(0) \in \mathbb R^n$, and let
$x(t)$ denote a maximal solution of the vector field. Since the
vector field has norm 1, $x(t)$ cannot diverge in finite time and
therefore $x(t)$ is well-defined for all $t \in \mathbb R$.
\par
Now we look at the function $t \mapsto y \cdot x(t) - f(x(t))$. Its
derivative w.r.t. $t$ is $(y - \nabla f_{x(t)}) \dot x(t) > \eps$.
Therefore, $\lim_{t \rightarrow \infty} y \cdot x(t) - f(x(t)) = \infty$.
We deduce from there that $\sup_{r \in \mathbb R^n} y \cdot r - f(r) = \infty$
Hence, $y \not \in L_f$, a contradiction.
\end{proof}

By replacing $f$ by $g_A$, we conclude that the image
of the momentum map $\nabla g_A$ is convex.

\begin{proof}[Proof of Proposition~\ref{convexity}]
The momentum map $\nabla g_A$ maps $\toric{n}$ onto the interior
of $\cvx A$. Indeed, let $a=A^{\alpha}$ be a row of $A$, associated to a
vertex of $\cvx A$. Then there is a direction $v \in \mathbb R^n$ 
such that
\[
a\cdot v = \max_{x \in \cvx A} x \cdot v
\]
for some unique $a$. 
\par
We claim that $a \in \overline{\nabla g_A (\mathbb R^n)}$. Indeed,
let $x(t) = v_A (tv)$, $t$ a real parameter. If $b$ is another 
row of $A$,
\[
e^{a\cdot tv} = e^ {t {a \cdot v}} \gg e^{t{b \cdot v}} = e^{b \cdot tv}
\]
as $t \rightarrow \infty$. We can then write $\hat v_A(tv)^{2T}$ as:
\[
\hat v_A (tv) = 
\left[
\begin{matrix}
\vdots \\
e^{t a \cdot v} \\
\vdots
\end{matrix}
\right]^T
C
\diag
\left[
\begin{matrix}
\vdots \\
e^{t a \cdot v} \\
\vdots
\end{matrix}
\right]
\punct{.}
\]
\par
Since $C$ is positive definite, $C_{\alpha \alpha}>0$ and
\[
\lim_{t \rightarrow \infty}
v_A (tv) ^{2T}
=
\lim_{t \rightarrow \infty}
\frac{ \hat v_A (tv)^{2T} }{\| \hat v_A (tv) \|^2} 
=
\mathrm{e}_a ^T \frac{C_{\alpha \alpha}}{C_{\alpha \alpha}}
=
\mathrm{e}_a ^T
\punct{,}
\]
where $\mathrm{e}_a$ is the unit vector in $\mathbb R^M$ corresponding to
the row $a$. 
It follows that $\lim_{t\rightarrow \infty} \nabla g_A(tv) = a$

When we set $q=0$, we have $\det D^2 g_A \ne 0$ on $\mathbb R^n$, so we
have a local diffeomorphism at each point $p \in \mathbb R^n$.
Assume that $(\nabla g_A)_p = (\nabla g_A)_{p'}$ for $p \ne p'$.
Then, let $\gamma (t) = (1-t)p + tp'$. The function 
$t \mapsto (\nabla g_A)_{\gamma(t)} \gamma'(t)$ has the same value at $0$
and at $1$, hence by Rolle's Theorem its derivative must vanish at some 
$t^* \in (0,1)$. 
\par
In that case,
\[
(D^2 g_A)_{\gamma(t^*)} (\gamma'(t^*), \gamma'(t^*)) = 0
\]
and since $\gamma'(t^*) = p'-p \ne 0$, $\det D^2g_A$ must vanish
in some $p \in \mathbb R^n$. This contradicts
Lemma~\ref{dimnondeg}.
\end{proof}

\subsection{More Properties of the Momentum Map}

\begin{figure}
\setlength{\unitlength}{2mm}
\centerline{
\scalebox{.5}{
\begin{picture}(50,54)(0,8)
\put(10,30){\line(0,1){20}}
\put(10,30){\circle*{0.5}}
\put(11,31){$A^1=(1,3)$}
\put(10,50){\line(1,1){10}}
\put(10,50){\circle*{0.5}}
\put(11,48){$A^2=(1,5)$}
\put(20,60){\line(1,-1){20}}
\put(20,60){\circle*{0.5}}
\put(21,61){$A^3=(2,6)$}
\put(40,40){\line(-1,-3){10}}
\put(40,40){\circle*{0.5}}
\put(41,40){$A^4=(4,4)$}
\put(30,10){\line(-1,1){20}}
\put(30,10){\circle*{0.5}}
\put(31,9){$A^5=(1,3)$}
\thicklines
\put(30,30){\vector(-1,3){10}}
\put(30,30){\circle*{0.5}}
\put(28,28){$\nabla g_A$}
\put(27,40){$A^3 - \nabla g_A$}
\end{picture}
}}
\caption{Geometric interpretation of $Dv_A^{\alpha}$\label{figmomentum}}
\end{figure}
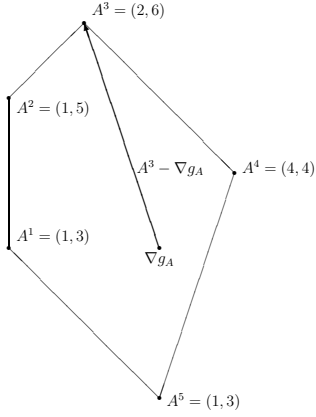

We can also give an interpretation of the derivative
$Dv_A$ in terms of the momentum map (see figure~\ref{figmomentum}).

\begin{lemma} \label{lem:geom}
\[
(Dv_A^{\alpha})_{p} u = |(v_A^{\alpha})_{p}| \left( A^{\alpha}- \nabla g_A(p) \right) \cdot u
\punct{.}
\]
\end{lemma}

where $|(v_A^{\alpha})_{p}|$ stands for $\frac{|(\hat v_A^{\alpha})_{p}|}{\|(\hat v_A)_{p}\|
}$,
and where $A^{\alpha}$ and $\nabla g_A$ are co-vectors.

\begin{proof}
By formula~\ref{deriv},
\[
(Dv_A)_{\scriptscriptstyle p} u =
\mathrm{diag} v_A(p) A u
-
v_A (p) 
v_A (p) ^T
\mathrm{diag} \left( v_A(p) \right) A u
\punct{.}
\]
Hence its $\alpha$-th coordinate is:
\[
(Dv_A^{\alpha})_{\scriptscriptstyle p} u =
(v_A^{\alpha})_{p}
\left(
A^{\alpha} u - (v_A)_{p}^{2T} A u
\right)
=
(v_A^{\alpha})_{p}
(A^{\alpha} - \nabla g_A(p)) u
\punct{.}
\]
\end{proof}

Bearing in mind that $\sum v_A^{\alpha}(p)^2 = 1$, we
obtain an immediate consequence:
\begin{lemma}\label{usefulbound}
For all $(p,q) \in \mathcal T^n$,
\[
\| Dv_A(p,q) \|
\le \mathrm{diam} (\cvx A)
\text{ \ and}
\]
\[
\| \frac{1}{2} D^2g_A(p) \|
\le \left( \mathrm{diam} (\cvx A) \right)^2
\]
\end{lemma}

\subsection{Evaluation Map and Condition Matrix}
\label{sec:bernshtein}
In the setting of Theorem~\ref{GEN-BERNSHTEIN}, we can 
identify each space of polynomials $(\mathcal F_{A_i},
\langle\cdot,\cdot\rangle_{\scriptscriptstyle A_i})$ to the (co)vector space
$(\mathbb C^{M_i})^{\vee}$, endowed with the canonical
inner product. The 
value of $f_i$ at $\exp(p+q\sqrt{-1})$ is then precisely 
$f_i \cdot \hat v_{A_i}(p+q\sqrt{-1})$.
\par
More generally, we can 
define the {\em evaluation map} by 
\[
\begin{array}{lrcl}
\mathit{ev}: & (\mathcal F_{A_1} \times \cdots \times \mathcal F_{A_n})
\times \mathcal T^n & \rightarrow & \mathbb C^n \\
   & (f_1, \cdots, f_n; p+q\sqrt{-1}) & \mapsto & 
\left[
\begin{matrix}
f^1 \cdot \hat v_{A_1}(p+q \sqrt{-1})\\ 
\vdots \\
f^n \cdot \hat v_{A_n}(p+q \sqrt{-1})
\end{matrix}
\right]
\end{array}
\punct{.}
\]
\medskip
\par

Following~\cite{BCSS}, we look at the linearization of the
implicit function $p+q\sqrt{-1} = G(f)$ for the equation
$\mathit{ev}(f,p+q\sqrt{-1})=0$.

\begin{definition} The {\em condition matrix} of $\mathit{ev}$ at
$(f, p+q\sqrt{-1})$ is
\[
DG = D_{\scriptscriptstyle \mathcal T^n}
(\mathit{ev})^{-1} D_{\scriptscriptstyle \mathcal F}
(\mathit{ev})
\punct{,}
\]
where $\mathcal F = \mathcal F_{A_1} \times \cdots \times \mathcal F_{A_n}$.
\end{definition}

Above, $D_{\scriptscriptstyle \mathcal T^n} (\mathit{ev})$ is a linear operator
from an $n$-dimensional complex space into $\mathbb C^n$, while 
$D_{\scriptscriptstyle \mathcal F} (\mathit{ev})$ goes 
from an $M_1+\cdots+M-n$-dimensional complex space into $\mathbb C^n$.

\begin{lemma} \label{herecomesthewedge}
  Assume that $\mathit{ev}(f;p+q\sqrt{-1}) = 0$. Then,
\begin{multline*}
\det \left( DG DG^H \right) ^{-1}
dp_1\wedge dq_1 \wedge \cdots \wedge dp_n\wedge dq_n 
=
(-1)^{n(n-1)/2}\ 
\bigwedge \sqrt{-1} 
f^i \cdot (Dv_{A_i})_{\scriptscriptstyle (p,q)}dp
\wedge \\
\wedge \bar f^i \cdot (Dv_{A_i})_{\scriptscriptstyle (p,-q)}dq
\punct{.}
\end{multline*}
\end{lemma}

Note that although $f^i \cdot (Dv_{A_i})_{\scriptscriptstyle (p,q)}dp$
is a complex-valued form, each wedge
$f^i \cdot (Dv_{A_i})_{\scriptscriptstyle (p,q)}dp
\wedge \bar f^i \cdot (Dv_{A_i})_{\scriptscriptstyle (p,-q)}dq$
is a real-valued 2--form.

\begin{proof}

We compute:
\[
D_{\scriptscriptstyle \mathcal F}(\mathit{ev})|_{(p,q)} 
= 
\left[
\begin{matrix}
\sum_{\alpha=1}^{M_1} \hat v_{A_1}^{\alpha} (p+q \sqrt{-1}) df^1_{\alpha}\\
\vdots \\
\sum_{\alpha=1}^{M_n} \hat v_{A_n}^{\alpha} (p+q \sqrt{-1}) df^n_{\alpha}
\end{matrix}
\right]
\punct{,}
\]

and hence
\[
D_{\scriptscriptstyle \mathcal F}
(\mathit{ev}) D_{\scriptscriptstyle \mathcal F}(\mathit{ev})^H = 
\mathrm{diag \ } \| \hat v_ {A_i} \|^2
\punct{.}
\]

Also,

\[
D_{\scriptscriptstyle \mathcal T^n}(\mathit{ev}) = 
\left[
\begin{matrix}
f^1 \cdot D\hat v_{A_1}\\
\vdots \\
f^n \cdot D\hat v_{A_n}
\end{matrix}
\right]
\punct{.}
\]

Therefore,

\[
\det \left(DG_{\scriptscriptstyle (p,q)} DG_{\scriptscriptstyle (p,q)}^H\right)^{-1}
=
\left| \det 
\left[
\begin{matrix}
f^1 \cdot \frac{1}{\|\hat v_{A_1}\|} D\hat v_{A_1} \\
\vdots \\
f^n \cdot \frac{1}{\|\hat v_{A_n}\|} D\hat v_{A_n} 
\end{matrix}
\right]
\right|^2
\punct{.}
\]

We can now use Lemma~\ref{lem-orthog} to conclude the following:

\formula{\label{det:cond:matrix}Determinant of the Condition Matrix}
{
\det \left(DG_{\scriptscriptstyle (p,q)} DG_{\scriptscriptstyle (p,q)}^H\right)^{-1}
=
\left| \det 
\left[
\begin{matrix}
f^1 \cdot Dv_{A_1}\\ 
\vdots \\
f^n \cdot Dv_{A_n} 
\end{matrix}
\right]
\right|^2
}

We can now write the same formula as a determinant of a block matrix:

\[
\det \left(DG_{\scriptscriptstyle (p,q)} DG_{\scriptscriptstyle (p,q)}^H\right)^{-1}
=
\det 
\left[
\begin{matrix}
f^1 \cdot Dv_{A_1}&\\ 
\vdots &\\
f^n \cdot Dv_{A_n}&\\ 
&\bar f^1 \cdot D\bar v_{A_1}\\ 
&\vdots \\
&\bar f^n \cdot D\bar v_{A_n} 
\end{matrix}
\right]
\]

and replace the determinant by a wedge.
The factor $(-1)^{n(n-1)/2}$ comes from replacing
$dp_1 \wedge \cdots \wedge dp_n \wedge dq_1 \wedge \cdots \wedge dq_n$
by $dp_1 \wedge dq_1 \wedge \cdots \wedge dp_n \wedge dq_n$.
\end{proof}

\begin{proof}[Proof of Theorem~\ref{GEN-BERNSHTEIN}]
Given $(p,q)\in \toric{n}$, we define $\mathcal F_{(p,q)}$ as
the space of $f \in \mathcal F_{A_1} \times \cdots \times 
\mathcal F_{A_n}$ such that $\mathit{ev}(f;p+q\sqrt{-1})$
vanishes.

Using \cite[Theorem 5 p.~243]{BCSS} (or Proposition~\ref{coarea2} 
p.~\pageref{coarea2} below), we deduce that the average number of 
complex roots is:
\[
\mathrm{Avg} =
\int_{(p,q) \in U}
\int_{f \in \mathcal F_{(p,q)}}
\left(\prod \frac{e^{-\|f_i\|^2/2}}{(2 \pi)^{M_i}}\right) 
\det \left(DG_{\scriptscriptstyle (p,q)} DG_{\scriptscriptstyle (p,q)}^H\right)^{-1}
\punct{.}
\]

By Lemma~\ref{herecomesthewedge}, we can replace the inner integral by
a $2n$--form valued integral:
\begin{multline*}
\mathrm{Avg} =
(-1)^{n(n-1)/2}\ 
\int_{(p,q) \in U}
\int_{f \in \mathcal F_{(p,q)}}
\bigwedge_i \frac{e^{-\|f_i\|^2/2}}{(2 \pi)^{M_i}}
f^i (Dv_{A_i})_{\scriptscriptstyle (p,q)}dp
\wedge \\
\wedge \bar f^i (Dv_{A_i})_{\scriptscriptstyle (p,-q)}dq
\punct{.}
\end{multline*}

Since the image of $Dv_{A_i}$ is precisely $\mathcal F_{A_i}|_{(p,q)} 
\subset \mathcal F_{A_i}$,
one can add $n$ extra variables corresponding to the directions 
$v_{A_i}(p+q \sqrt{-1})$
without changing the integral: we write $\mathcal F_{A_i} = \mathcal 
F_{A_i,(p,q)} \times
\mathbb C v_{A_i}(p+q \sqrt{-1})$. Since $\left( f^i + t v_{A_i}(p+q \sqrt{-1}) \right)
Dv_{A_i}$ is equal to $f^i Dv_{A_i}$, the average number of roots is indeed:
\begin{multline*}
\mathrm{Avg} =
(-1)^{n(n-1)/2}\ 
\int_{(p,q) \in U}
\int_{f \in \mathcal F}
\bigwedge_i \frac{e^{-\|f_i\|^2/2}}{(2 \pi)^{M_i +1}}
f^i \cdot (Dv_{A_i})_{\scriptscriptstyle (p,q)}dp
\wedge \\
\wedge \bar f^i \cdot (Dv_{A_i})_{\scriptscriptstyle (p,-q)}dq
\punct{.}
\end{multline*}

In the integral above, all the terms that are multiple of $f^i_{\alpha} \bar f^i_{\beta}$
for some $\alpha \ne \beta$ will cancel out. Therefore,
\begin{multline*}
\mathrm{Avg} =
(-1)^{n(n-1)/2}\ 
\int_{(p,q) \in U}
\int_{f \in \mathcal F}
\bigwedge_i \frac{e^{-\|f_i\|^2/2}}{(2 \pi)^{M_i+1}}
\sum_{\alpha} |f^i_{\alpha}|^2 (Dv_{A_i})^{\alpha}_{\scriptscriptstyle (p,q)}dp
\wedge 
\\
\wedge
(Dv_{A_i})^{\alpha}_{\scriptscriptstyle (p,-q)}dq
\punct{.}
\end{multline*}

Now, we apply the integral formula:
\[
\int_{x \in \mathbb C^{M} } |x_1|^2 \frac{e^{-\|x\|^2 / 2}}{(2 \pi)^{M}}
=
\int_{x_1 \in \mathbb C} |x_1|^2 \frac{e^{-|x_1|^2 / 2}}{2 \pi}
= 2
\]
to obtain:
\[
\mathrm{Avg} =
\frac{(-1)^{n(n-1)/2}}{\pi^n}
\int_{(p,q) \in U}
\bigwedge \sum_{\alpha} (Dv_{A_i})^{\alpha}_{\scriptscriptstyle (p,q)}dp
\wedge (Dv_{A_i})^{\alpha}_{\scriptscriptstyle (p,-q)}dq
\punct{.}
\]

According to formul\ae~\ref{secder} and~\ref{omega}, the integrand is just 
$2^{-n} \bigwedge \omega_{A_i}$, and thus
\[
\mathrm{Avg} =
\frac{(-1)^{n(n-1)/2}}{\pi^n}
\int_{U}
\bigwedge_i \omega_{A_i}
=
\frac{n!}{\pi^n}
\int_{U}
d\toric{n}
\punct{.}
\]
\end{proof}

\section{The Condition Number}

\subsection{Proof of Theorem~\ref{condnumber}}

Let $(p,q) \in \toric{n}$ and let
$f \in \mathcal F_{(p,q)}$. Without loss of generality, we can
assume that $f$ is scaled so that for all $i$,
$\|f^i\|=1$. 
\par
Let $\delta f \in \mathcal F_{(p,q)}$
be such that $f+\delta f$ is singular at $(p,q)$, and assume
that $\sum \| \delta f^i \|^2$ is minimal. Then, due
to the scaling we chose,
\[
d_{\mathbb P} (f, \Sigma_{(p,q)}) = \sqrt{\sum \| \delta f^i \|^2} 
\punct{.}
\]
\par
Since $f + \delta f$ is singular, there is a vector $u \ne 0$ such
that
\[
\left[
\begin{matrix}
(f^1 + \delta f^1) \cdot (D\hat v_{A_1})_{(p,q)}
\\
\vdots
\\
(f^n +\delta f^n) \cdot (D\hat v_{A_n})_{(p,q)}
\end{matrix}
\right]
u = 0
\]
and hence
\[
\left[
\begin{matrix}
(f^1 + \delta f^1) \cdot (Dv_{A_1})_{(p,q)}
\\
\vdots
\\
(f^n +\delta f^n) \cdot (Dv_{A_n})_{(p,q)}
\end{matrix}
\right]
u = 0
\punct{.}
\]
\medskip
\par 
This means that
\[
\left\{
\begin{array}{lcl}
f^1 \cdot Dv_{A_1} u &=& - \delta f^1 \cdot Dv_{A_1} u \\
&\vdots&\\
f^n \cdot Dv_{A_n} u &=& - \delta f^n \cdot Dv_{A_n} u \\
\end{array}
\right.
\punct{.}
\]
\par
Let $D(f)$ denote the matrix
\[
D(f) \defeq
\left[
\begin{matrix}
f^1 \cdot (Dv_{A_1})_{(p,q)}
\\
\vdots
\\
f^n \cdot (Dv_{A_n})_{(p,q)}
\end{matrix}
\right]
\punct{.}
\]

Given $v = D(f) \ u$, we obtain:
\begin{equation}\label{eqonv}
\left\{
\begin{array}{lcl}
v_1 &=& - \delta f^1 \cdot Dv_{A_1} D(f)^{-1} v \\
&\vdots&\\
v_n &=& - \delta f^n \cdot Dv_{A_n} D(f)^{-1} v \\
\end{array}
\right.
\end{equation}
\par

We can then scale $u$ and $v$, such that $\|v\|=1$.

\begin{claim}
Under the assumptions above, $\delta f^i$ is colinear to 
$\left( Dv_{A_i} D(f)^{-1} v \right)^H$.
\end{claim}

\begin{proof}
Assume that $\delta f^i = g + h$, with $g$ colinear and $h$ orthogonal to
$\left( Dv_{A_i} D(f)^{-1} v \right)^H$.
As the image of $Dv_{A_i}$ is orthogonal to $v_{A_i}$, $g$ is orthogonal
to $v_{A_i}^H$, so $\mathit{ev}(g^i,(p,q))=0$ and hence
$\mathit{ev}(h^i,(p,q))=0$. We can therefore
replace $\delta f^i$ by $g$ without compromising
equality~(\ref{eqonv}). Since $\|\delta f\|$ was minimal, this implies
$h=0$.
\end{proof}

We obtain now an explicit expression for $\delta f^i$ in terms of $v$:

\[
\delta f^i
= - v_i
\frac{
 \left( Dv_{A_i} D(f)^{-1} v \right)^H
}
{
\| Dv_{A_i} D(f)^{-1} v \|^2
}
\punct{.}
\]

Therefore,
\[
\| \delta f^i \|
=
\frac{ |v_i| }
{\|Dv_{A_i} D(f)^{-1} v \|}
=
\frac{ |v_i| }
{\|\left( D(f)^{-1} v \right)\|_{A_i}}
\punct{.}
\]

So we have proved the following result: 

\begin{lemma} \label{lemonv}
Fix $v$ so that $\|v\|=1$ and 
let $\delta f \in \mathcal F_{(p,q)}$ be such that equation (\ref{eqonv}) 
holds and $\|\delta f\|$ is minimal. Then, 
\[
\| \delta f^i \|
=
\frac{ |v_i| }
{\|D(f)^{-1} v \|_{A_i}}
\punct{.}
\]
\end{lemma}

Lemma~\ref{lemonv} provides an immediate lower bound for 
$\| \delta f \| = \sqrt{\sum \| \delta f^i\|^2}$: Since
\[
\| \delta f^i \| \ge
\frac{ |v_i| }
{\max_j \| D(f)^{-1} v \|_{A_j}}
\punct{,}
\]
we can use $\|v\|=1$ to deduce that
\[
\sqrt{\sum_i \| \delta f^i \|^2} \ge
\frac{ 1 }
{\max_j \|D(f)^{-1} v \|_{A_j}}
\ge
\frac{ 1 }
{\max_j \| D(f)^{-1} \|_{A_j}}
\punct{.}
\]

Also, for any $v$ with $\|v\|=1$, we can choose $\delta f$ minimal so
that equation~(\ref{eqonv}) applies. Using Lemma~\ref{lemonv}, we
obtain:
\[
\| \delta f^i \| \le
\frac{ |v_i| }
{\min_j \| D(f)^{-1} v \|_{A_j}}
\punct{.}
\]

Hence
\[
\sqrt{\sum_i \| \delta f^i \|^2} \le
\frac{ 1 }
{\min_j \| D(f)^{-1} v \|_{A_j}}
\punct{.}
\]

Since this is true for any $v$, and $\| \delta f\|$ is
minimal for all $v$, we have 
\[
\sqrt{\sum_i \| \delta f^i \|^2} \le
\frac{ 1 }
{\max_{\|v\|=1} \min_j \| D(f)^{-1}\|_{A_j}}
\]

and this proves Theorem~\ref{condnumber}.

\subsection{Idea of the Proof of Theorem~\ref{mixed-cond}}

The proof of Theorem~\ref{mixed-cond} is long. 
We first sketch the idea of the proof. Recall
that $\mathcal F_{(p,q)}$ is the set of all
$f \in \mathcal F$ such that $\mathit{ev}(f;p
+q\sqrt{-1}) = 0$, and that $\Sigma_{(p,q)}$
is the restriction of the discriminant to the
fiber $\mathcal F_{(p,q)}$:
\[
\Sigma_{(p,q)} \defeq \{ f \in \mathcal F_{(p,q)}:
 D(f)_{(p,q)} \text{\ does not have full rank} \} 
\punct{.}
\]
\medskip
\par
The space $\mathcal F$ is endowed with a Gaussian
probability measure, with volume element 
\[
\frac{e^{-\|f\|^2 / 2}}{(2 \pi)^{\sum M_i}} d\mathcal F
\punct{,}
\]
where $d\mathcal F$ is the usual volume form in 
$\mathcal F = (\mathcal F_{A_1}, \langle\cdot,\cdot
\rangle_{\scriptscriptstyle A_1}) \times \cdots
\times (\mathcal F_{A_n}, \langle\cdot,\cdot\rangle_{\scriptscriptstyle A_n})$
and $\|f\|^2 = \sum \|f^i\|_{A_i}^2$.
For $U$ a set in $\toric{n}$, we defined earlier (in the statement of 
Theorem \ref{mixed-cond}) the quantity: 
\[
\nu^{A} (U,\eps) \defeq 
\mathrm{Prob} [ \mymu(f,U) > \eps^{-1}]
=
\mathrm{Prob} [ 
\exists (p,q) \in U: d_{\mathbb P}(f,\Sigma_{(p,q)}) < \eps]
\punct{.}
\]

The na\"{\i}ve idea for bounding $\nu^A(U,\eps)$ is as follows:
Let $V(\eps) \defeq \{ (f,(p,q)) \in \mathcal F \times U: 
\mathit{ev}(f;(p,q))=0 \text{\ and \ }
d_{\mathbb P} (f, \Sigma_{(p,q)}) < \eps \}$. We also define $\pi: V(\eps)
\rightarrow \mathcal F$ as the canonical projection mapping 
$\cF\times U$ to $\cF$, and set 
$\#_{V(\eps)}(f) \defeq \# \{ (p,q) \in U: (f,(p,q)) \in V(\eps) \}$.
Then,

\begin{eqnarray*}
\nu^A(U,\eps)
&=&
\int_{f \in \mathcal F}
\chi_{\pi(V(\eps))} (f) \ 
\frac{e^{-\|f\|^2 / 2}}{(2 \pi)^{\sum M_i}} d\mathcal F
\\
&\le&
\int_{f \in \mathcal F}
\#_{V(\eps)} \ 
\frac{e^{-\|f\|^2 / 2}}{(2 \pi)^{\sum M_i}} d\mathcal F
\end{eqnarray*}
with equality in the linear case.
\par
Now we apply the coarea formula~\cite[Theorem 5 p.~243]{BCSS} to obtain:

\[
\nu^A(U,\eps)
\le
\int_{(p,q) \in U \subset \toric{n}}
\int_{
\substack{
f \in \mathcal F_{(p,q)} \\
d_{\mathbb P} (f, \Sigma_{(p,q)}) < \eps
}
}
\frac{1}{NJ(f;(p,q))} \ 
\frac{e^{-\|f\|^2 / 2}}{(2 \pi)^{\sum M_i}} d\mathcal F
\ dV_{\toric{n}}
\punct{,}
\]

where $dV_{\toric{n}}$ stands for Lebesgue measure in
$\toric{n}$. Again, in the linear case, we have
equality.

We already know from Lemma~\ref{herecomesthewedge}
that 
\[
1/NJ(;(p,q))
=
\bigwedge_{i=1}^n 
f^i \cdot (Dv_{A_i})_{(p,q)} dp \wedge \bar f^i \cdot (D\bar v_{A_i})_{(p,q)} dq
\punct{.}
\]
\par

We should focus now on the inner integral. In each coordinate space
$\mathcal F_{A_i}$, we can introduce a new orthonormal system of coordinates
(depending on $(p,q)$) by decomposing:
\[
f^i = f^i_{\romone} + f^i_{\romtwo} + f^i_{\romthree}
\punct{,}
\]

where $f^i_{\romone}$ is the component colinear
to $v_{A_i}^H$, 
$f^i_{\romtwo}$ is the projection of
$f^i$ to $({\mathrm{range}}\  Dv_{A_i})^H$,
and
$f^i_{\romthree}$ is orthogonal to 
$f^i_{\romone}$ and $f^i_{\romtwo}$.

Of course, $f^i \in (\mathcal F_{A_i})_{(p,q)}$ if and only if
$f^i_{\romone} = 0$.

Also,
\begin{multline*}
\bigwedge_{i=1}^n
f^i \cdot (Dv_{A_i})_{(p,q)} dp \wedge \bar f^i \cdot (D\bar v_{A_i})_{(p,q)} dq
=
\\
=
\bigwedge_{i=1}^n
f_{\romtwo}^i \cdot (Dv_{A_i})_{(p,q)} dp \wedge \bar f_{\romtwo}^i \cdot 
(D\bar v_{A_i})_{(p,q)} dq \punct{.}
\end{multline*}

It is an elementary fact that 
\[
d_{\mathbb P}(f^i_{\romtwo}+f^i_{\romthree}, \Sigma_{(p,q)})
\le
d_{\mathbb P}(f^i_{\romtwo}, \Sigma_{(p,q)})
\punct{.}
\]

It follows that for $f \in \mathcal F_{(p,q)}$:
\[
d_{\mathbb P}(f, \Sigma_{(p,q)})
\le
d_{\mathbb P}(f_{\romtwo}, \Sigma_{(p,q)})
\punct{,}
\]
with equality in the linear case. Hence, we obtain:

\begin{multline*}
\nu^A(U,\eps)
\le
\int_{(p,q) \in U \subset \toric{n}}
\int_{
\substack{
f \in \mathcal F_{(p,q)} \\
d_{\mathbb P} (f_{\romtwo}, \Sigma_{(p,q)}) < \eps
}}
\left(
\bigwedge_{i=1}^n
f_{\romtwo}^i \cdot (Dv_{A_i})_{(p,q)} dp \wedge \bar f_{\romtwo}^i \cdot (D\bar v_{A_i})_{(p,q)} dq
\right)
\cdot
\\
\cdot
\frac{e^{-\|f^i_{\romtwo}+f^i_{\romthree}\|^2 / 2}}{(2 \pi)^{\sum M_i}} d\mathcal F
\ dV_{\toric{n}}
\punct{,}
\end{multline*}
with equality in the linear case. We can integrate the 
$\sum (M_i - n - 1)$ variables $f_{\romthree}$
to obtain:

\begin{proposition}\label{assimpleasitgets}
\begin{multline*}
\nu^A(U,\eps)
\le
\int_{(p,q) \in U \subset \toric{n}}
\int_{
\substack
{
f_{\romtwo} \in \mathbb C^{n^2} \\
d_{\mathbb P} (f_{\romtwo}, \Sigma_{(p,q)}) < \eps
}}
\left(
\bigwedge_{i=1}^n
f_{\romtwo}^i \cdot (Dv_{A_i})_{(p,q)} dp \wedge \bar f_{\romtwo}^i \cdot (D\bar v_{A_i})_{(p,q)} dq
\right)
\cdot \\ \cdot
\frac{e^{-\|f^i_{\romtwo}\|^2 / 2}}{(2 \pi)^{n(n+1)}} 
\ dV_{\toric{n}}
\punct{.}
\end{multline*}
with equality in the linear case. \qed 
\end{proposition}

\subsection{From Gaussians to Multiprojective Spaces}

The domain of integration in Proposition~\ref{assimpleasitgets}
makes integration extremely difficult. In order to estimate
the inner integral, we will need to perform a change of
coordinates.

Unfortunately, the Gaussian in Proposition~\ref{assimpleasitgets}
makes that change of coordinates extremely hard, and we will have
to restate Proposition \ref{assimpleasitgets} in terms of integrals over 
a product of projective spaces.  

The domain of integration will be $\mathbb P^{n-1} \times 
\cdots \times \mathbb P^{n-1}$. Translating an integral
in terms of Gaussians to an integral in terms of
projective spaces is not immediate, and we will use
the following elementary fact about Gaussians:

\begin{lemma}\label{gaussians1}
 Let $\varphi: \mathbb C^n \rightarrow \mathbb R$ be
$\mathbb C^*$-invariant (in the sense of the usual 
scaling action). Then we can also interpret
$\varphi$ as a function from $\mathbb P^{n-1}$ into
$\mathbb R$, and:
\[
\frac{1}{\vol(\mathbb P^{n+1})} 
\int_{[x]\in \mathbb P^{n-1}} \varphi(x)d[x]  =
\int_{x \in \mathbb C^{n}}  \varphi(x) 
\frac{e^{-\|x\|^2/2}}{(2\pi)^n} dx 
\punct{,}
\]
where, respectively, the natural volume forms on $\Pro^{n-1}$ and $\Cn$ 
are understood for each integral. 
\end{lemma}

Now the integrand in Proposition~\ref{assimpleasitgets}
is not $\mathbb C^*$--invariant. This is why we will need 
the following formula: 

\begin{lemma}\label{gaussians2}
 Under the hypotheses of Lemma~\ref{gaussians1},
\[
\frac{1}{\vol(\mathbb P^{n+1})} 
\int_{[x]\in \mathbb P^{n-1}} \varphi(x) d[x] =
\frac{1}{2n}
\int_{x \in \mathbb C^{n}}  \|x\|^2 \varphi(x) 
\frac{e^{-\|x\|^2/2}}{(2\pi)^n} dx
\punct{.}
\]
where, respectively, the natural volume forms on $\Pro^{n-1}$ and $\Cn$ 
are understood for each integral. 
\end{lemma}

\begin{proof}
\begin{eqnarray*}
\int_{x \in \mathbb C^{n}}  \|x\|^2 \varphi(x) 
\frac{e^{-\|x\|^2/2}}{(2\pi)^n} dx  
&=&
\int_{\Theta \in S^{2n-1}} 
\int_{r=0}^{\infty}
|r|^{2n+1} \varphi(\Theta) 
\frac{e^{-|r|^2/2}}{(2\pi)^n} dr d\Theta 
\end{eqnarray*} 
\[=
\int_{\Theta \in S^{2n-1}} 
\left(
-
\left[
|r|^{2n} 
\frac{e^{-|r|^2/2}}{(2\pi)^n} 
\right]_0^{\infty}
\right. 
\left. 
+
2n
\int_{r=0}^{\infty}
|r|^{2n-1} 
\frac{e^{-|r|^2/2}}{(2\pi)^n} 
dr 
\right)
\varphi(\Theta) 
d\Theta \] 
\[=
2n
\int_{x \in \mathbb C^{n}} \varphi(x) 
\frac{e^{-\|x\|^2/2}}{(2\pi)^n} dx 
\] 
\end{proof}

We can now introduce the notation:

\[
\mathrm{WEDGE}^A(f_{\romtwo}) 
\defeq 
\bigwedge_{i=1}^n 
\frac{1}{\|f_{\romtwo}^i\|^2} 
f_{\romtwo}^i \cdot (Dv_{A_i})_{(p,q)} dp \wedge 
\bar f_{\romtwo}^i \cdot (D\bar v_{A_i})_{(p,q)} dq
\punct{.}
\]
\par
This function is invariant under the $(\mathbb C^*)^n$-action
$\lambda \star f_{\romtwo} : f_{\romtwo} \mapsto 
(\lambda_1 f_{\romtwo}^1, \cdots, \lambda_n f_{\romtwo}^n)$. 

We adopt the following conventions:
$\mathcal F_{\romtwo} \subset \mathcal F$ is the space spanned by
coordinates $f_{\romtwo}$ and $\mathbb P (\mathcal F_{\romtwo}) $
is its quotient by $(\mathbb C^*)^n$.

We apply $n$ times
Lemma~\ref{gaussians2} and obtain:

\begin{proposition}\label{multiprojective}
Let $\mathrm{VOL} \defeq \vol(\mathbb P^{n-1})^{n}$. Then,
\[
\nu^A(U,\eps)
\le
\frac{(2n)^n}{\mathrm{VOL}}
\int_{(p,q) \in U \subset \toric{n}}
\int_{
\substack{
f_{\romtwo} \in \mathbb P (\mathcal F_{\romtwo})\\
d_{\mathbb P} (f_{\romtwo}, \Sigma_{(p,q)}) < \eps
} }
\mathrm{WEDGE}^A(f_{\romtwo}) \
d\mathbb P (\mathcal F_{\romtwo})
\ dV_{\toric{n}}
\]
and in the linear case,
\[
\nu^{\linear}(U,\eps)
=
\frac{(2n)^n}{\mathrm{VOL}}
\int_{(p,q) \in U \subset \toric{n}}
\int_{
\substack{
g_{\romtwo} \in \mathbb P (\mathcal F_{\romtwo}^{\linear}) \\
d_{\mathbb P} (g_{\romtwo}, \Sigma_{(p,q)}^{\linear}) < \eps
} }
\mathrm{WEDGE}^{\linear}(g_{\romtwo}) \
d\mathbb (P \mathcal F_{\romtwo}^{\linear})
dV_{\toric{n}} \ \text{\qed} 
\]
\end{proposition}

Now we introduce the following change of coordinates. Let
$L \in GL(n)$ be such that the minimum in Definition~\ref{intrinsic}
p.~\pageref{intrinsic} is attained:
\[
\begin{array}{crcl}
 \varphi: &
\mathbb P^{n-1} \times \cdots \times \mathbb P^{n-1} &
\rightarrow &
\mathbb P^{n-1} \times \cdots \times \mathbb P^{n-1} \\
&
f_{\romtwo} & \mapsto & g_{\romtwo} \defeq \varphi(f_{\romtwo}) \punct{,}\text{\ such}\\
& & & \text{that\ } 
g_{\romtwo}^i = f_{\romtwo}^i \cdot Dv_{A_i} L \punct{.}
\end{array}
\]

Without loss of generality, we scale $L$ such that $\det L=1$.
The following property follows from the definition of
$\mathrm{WEDGE}$:

\begin{equation}\label{wedge1}
\mathrm{WEDGE}^A (f_{\romtwo})
=
\mathrm{WEDGE}^{\linear} (g_{\romtwo}) \ 
\prod_{i=1}^n
\frac {\| g_{\romtwo}^i \|^2} {\| f_{\romtwo}^i \|^2}
\end{equation}

Assume now that
$d_{\mathbb P} (f_{\romtwo}, \Sigma_{(p,q)}) < \eps$.
Then there is $\delta f \in \mathcal F_{\romtwo}$, 
such that $f + \delta f \in \Sigma_{(p,q)}^{\linear}$ and
$\| \delta f \| \le \eps$ (assuming the scaling  $\|f^i_{\romtwo}\|=1$
for all $i$).
\par
Setting $g_{\romtwo} = \varphi (f_{\romtwo})$ and $\delta g = \varphi(g)$,
we obtain that $g + \delta g \in \Sigma_{(p,q)}^{\linear}$. 
\[
d_{\mathbb P} (g, \Sigma_{(p,q)}^{\linear})
\le
\sqrt{ \sum_{i=1}^n \frac {\| \delta g^i \|^2}{\| g^i_{\romtwo} \|^2}
}
\]

At each value of $i$,
\[
\frac {\| \delta g^i \|}{\| g^i_{\romtwo} \|}
\le
\frac {\| \delta f^i \|}{\| f^i_{\romtwo} \|}
\kappa ( D_{f^i_{\romtwo}}\varphi^i) 
\]
where $\kappa$ denotes Wilkinson's condition number of the
linear operator $D_{f^i_{\romtwo}}\varphi^i$. This is precisely
$\kappa (Dv_{A_i} L)$. Thus,
\[
d_{\mathbb P} (g, \Sigma_{(p,q)}^{\linear})
\le
\eps
\max_i \kappa (Dv_{A_i} L) = \max_i \sqrt{\kappa ( \omega_{A_i} )}
\]
 
Thus, an $\eps$-neighborhood of $\Sigma^A_{(p,q)}$ is mapped
into a $\sqrt{\kappa_U} \eps$ neighborhood of
$\Sigma^{\linear}_{(p,q)}$.

We use this property and equation~(\ref{wedge1}) to bound:

\begin{multline}\label{wedge2}
\nu^A(U,\eps)
\le
\frac{(2n)^n}{\mathrm{VOL}}
\int_{(p,q) \in U \subset \toric{n}}
\int_{
\substack{
g_{\romtwo} \in \mathbb P^{n-1} \times \cdots \times \mathbb P^{n-1} \\
d_{\mathbb P} (g_{\romtwo}, \Sigma^{\linear}_{(p,q)}) < \sqrt{\kappa_U} \eps
} }
\mathrm{WEDGE}^{\linear}(g_{\romtwo})
\cdot
\\
\cdot
\prod_{i=1}^n
\frac {\| g_{\romtwo}^i \|^2} {\| f_{\romtwo}^i \|^2}
|J_{g_{\romtwo}} \varphi^{-1}|^2
\
d(\mathbb P^{n-1}\times \cdots \times \mathbb P^{n-1})
\ dV_{\toric{n}}
\end{multline}

where $J_{g_{\romtwo}} \varphi^{-1}$ is the Jacobian of $\varphi^{-1}$ at
$g_{\romtwo}$. 
\begin{remark}
Considering each $Dv_{A_i}$ as a map from $\mathbb C^n$
into $\mathbb C^n$, the Jacobian is:
\[
J_{g_{\romtwo}} \varphi^{-1}
=
\prod_{i=1}^n \frac{\|\varphi^{-1}(g_{\romtwo})^i\|^n}{\|g_{\romtwo}^i\|^n} 
\left( \det Dv_{A_i}^H Dv_{A_i} \right)^{-1/2}
\punct{.}
\]
We will not use this value in the sequel. \qed 
\end{remark}

In order to simplify the expressions for the bound on $\nu^A (U,\eps)$, it
is convenient to introduce the following notations:

\begin{eqnarray*}
dP
&\defeq&
\frac{(2n)^n}{\mathrm {VOL}}
\mathrm{WEDGE}^{\linear}(g_{\romtwo})
\frac{ d (\mathbb P^{n-1}\times \cdots \times \mathbb P^{n-1}) }
{n! \ (\omega_{\linear})^{\bigwedge n}}
\\
H &\defeq&
\prod_{i=1}^n
\frac {\| g_{\romtwo}^i \|^2} {\| f_{\romtwo}^i \|^2}
|J_g \varphi^{-1}|^2
\\
\chi_{\delta} &\defeq& \chi_{ \left\{ g: d_{\mathbb P}(g, \Sigma_{(p,q)}^{\linear}) < \delta \right\} }
\end{eqnarray*}

Now equation~(\ref{wedge2}) becomes:

\begin{equation}\label{wedge3}
\nu^A(U,\eps)
\le
n!
\int_{(p,q) \in U \subset \toric{n}}
(\omega_{\linear})^{\bigwedge n}
\int_{ g_{\romtwo} \in \mathbb P^{n-1} \times \cdots \times \mathbb P^{n-1} }
dP \
H(g_{\romtwo})
\
\chi_{\sqrt{\kappa_U} \eps} (g_{\romtwo})
\end{equation}

\begin{lemma}
 Let $(p,q)$ be fixed. Then $\mathbb P^{n-1} \times \cdots \times \mathbb P^{n-1}$
together with density  function $dP$, is a probability space.
\end{lemma}

\begin{proof}
  The expected number of roots in $U$ for a linear system is 
\[
n!
\int_{(p,q) \in U}
\omega_{\linear}^{\bigwedge n}
\int_{ g_{\romtwo} \in \mathbb P^{n-1} \times \cdots \times \mathbb P^{n-1} }
dP
\punct{.}
\]

It is also $n! \int_U \omega_{\linear}^{\bigwedge n}$. This holds for all $U$,
hence the volume forms are the same and
\[
\int_{ g_{\romtwo} \in \mathbb P^{n-1} \times \cdots \times \mathbb P^{n-1} }
dP
= 1
\punct{.}
\]
\end{proof}

This allows us to interpret the inner integral of equation~(\ref{wedge3})
as the expected value of a product. This is less than the product of
the expected values, and:

\begin{multline*}
\nu^A(U,\eps)
\le
n!
\int_{(p,q) \in U \subset \toric{n}}
(\omega_{\linear})^{\bigwedge n}
\left(
\int_{ g_{\romtwo} \in \mathbb P^{n-1} \times \cdots \times \mathbb P^{n-1} }
dP \
H(g_{\romtwo})
\right)
\cdot
\\
\cdot
\left(
\int_{ g_{\romtwo} \in \mathbb P^{n-1} \times \cdots \times \mathbb P^{n-1} }
dP \
\chi_{\sqrt{\kappa_U} \eps} (g_{\romtwo})
\right)
\end{multline*}

Because generic systems of linear equations have one root, we can also
consider $U$ as a probability space, with probability measure
$\frac{1}{\vol^{\linear} U} n! \omega_{\linear}^{\bigwedge n}$.
Therefore, we can bound:
\begin{multline*}
\nu^A(U,\eps)
\le
\frac{1}{\vol^{\linear} U}
\left(
\int_{(p,q) \in U}
n!
(\omega_{\linear})^{\bigwedge n}
\int_{ g_{\romtwo} \in \mathbb P^{n-1} \times \cdots \times \mathbb P^{n-1} }
dP \
H(g_{\romtwo})
\right)
\cdot \\ \cdot
\left(
\int_{(p,q) \in U}
n!
(\omega_{\linear})^{\bigwedge n}
\int_{ g_{\romtwo} \in \mathbb P^{n-1} \times \cdots \times \mathbb P^{n-1} }
dP \
\chi_{\sqrt{\kappa_U} \eps} (g_{\romtwo})
\right)
\end{multline*}

The first parenthesis is $\vol^A (U)$. The second parenthesis is $\nu^{\linear}(\sqrt{\kappa_U}\eps,U)$.
This concludes the proof of Theorem~\ref{mixed-cond}.

\begin{proof}[Proof of Corollary~\ref{unmixed:complex}]
We set $L = Dv_{A_i}^{\dagger}|{\mathrm{range} Dv_{A_i}}$,
then $\kappa (\omega_{A_1}, \cdots, \omega_{A_n}; (p,q)) = 1$.
\end{proof}

\section{Real Polynomials}

\subsection{Proof of Theorem~\ref{expected-real}}

\begin{proof}[Proof of Theorem~\ref{expected-real}]

As in the complex case (Theorem~\ref{GEN-BERNSHTEIN}),
the expected number of roots can be computed by
applying the co-area formula:
\[
AVG =
\int_{p \in U}
\int_{f \in \mathcal F_p^{\mathbb R}}
\prod_{i=1}^n 
\frac{e ^{-\|f^i\|^2/2}}{\sqrt{2 \pi}^{M_i}}
\det (DG \ DG^H)^{-1/2}
\punct{.}
\]

Now there are three big diferences. The set $U$ is
in $\mathbb R^n$ instead of $\toric{n}$, the
space $\mathcal F_p^{\mathbb R}$ contains
only real polynomials (and therefore has half the
dimension), and we are integrating the square root
of $1/ \det (DG\ DG^H)$.

Since we do not know in general how to integrate
such a square root, we bound the inner integral
as follows. We consider the real Hilbert space of
functions integrable in $\mathcal F_p^{\mathbb R}$
endowed with Gaussian probability measure. The
inner product in this space is:

\[
\langle \varphi , \psi\rangle \defeq
\int_{\mathcal F_p^{\mathbb R}}
\varphi(f) \psi(f)
\prod_{i=1}^n 
\frac{e^{-\|f^i\|^2/2}}{\sqrt{2 \pi}^{M_i-1}}
dV
\punct{,}
\]

where $dV$ is Lebesgue volume. If $\mathbf 1$ denotes
the constant function equal to $1$, we interpret

\[
AVG
=
\int_{p \in U}
(2\pi)^{-n/2}
\left\langle
\det (DG \ DG^H)^{-1/2}
,
\mathbf 1
\right\rangle
\punct{.}
\]

Hence Cauchy-Schwartz inequality implies:

\[
AVG
\le
\int_{p \in U}
(2\pi)^{-n/2}
\|
\det (DG \ DG^H)^{-1/2}
\|
\|
\mathbf 1
\|
\punct{.}
\]

By construction, $\| \mathbf 1\|=1$, and we are left with:

\[
AVG
\le
\int_{p \in U}
(2\pi)^{-n/2}
\sqrt{
\int_{\mathcal F_p^{\mathbb R}}
\prod_{i=1}^n 
\frac{e^{-\|f^i\|^2/2}}{\sqrt{2 \pi}^{M_i-1}}
\det (DG \ DG^H)^{-1}
}
\punct{.}
\]

As in the complex case, we add extra $n$ variables:
\[
AVG
\le
(2\pi)^{-n/2}
\int_{p \in U}
\sqrt{
\int_{\mathcal F^{\mathbb R}}
\prod_{i=1}^n 
\frac{e^{-\|f^i\|^2/2}}{\sqrt{2 \pi}^{M_i}}
\det (DG \ DG^H)^{-1}
}
\punct{,}
\]

and we interpret $\det (DG \ DG^H)^{-1}$ in terms of a wedge.
Since
\[
\int_{x \in \mathbb R^M} 
|x_1|^2 \frac{e^{-\|x\|^2/2}}{\sqrt{2 \pi}^M}
=
\int_{y \in \mathbb R} 
y^2 \frac{e^{-y^2/2}}{\sqrt{2 \pi}}
=
\int_{y \in \mathbb R} 
\frac{e^{-y^2/2}}{\sqrt{2 \pi}}
=1
\punct{,}
\]
we obtain:
\[
AVG
\le
(2\pi)^{-n/2}
\int_{p \in U}
\sqrt{
n! d\toric{n}
}
=
(2\pi)^{-n/2}
\int_{p \in U}
\sqrt{
n! d\toric{n}
}
\punct{.}
\]

Now we would like to use Cauchy-Schwartz again. This time,
the inner product is defined as:
\[
\langle \varphi , \psi \rangle \defeq \int_{p \in U} \varphi(p) \psi(p) dV
\punct{.}
\]

Hence,
\[
AVG
\le
(2\pi)^{-n/2}
\langle n! d\toric{n} , \mathbf{1} \rangle
\le
(2\pi)^{-n/2}
\| n! d\toric{n} \| \|\mathbf{1}\|
\punct{.}
\]

This time, $\|\mathbf{1}\|^2 = \lambda(U)$, so
we bound:

\begin{eqnarray*}
AVG
&\le&
(2\pi)^{-n/2}
\sqrt{\lambda(U)}
\sqrt{
\int_U n! d\toric{n}
}
\\
&\le&
(4\pi^2)^{-n/2}
\sqrt{\lambda(U)}
\sqrt{
\int_{(p,q)\in \toric{n}, p\in U} n! d\toric{n}
}
\punct{.}
\end{eqnarray*}

\end{proof}
\subsection{Proof of Theorem~\ref{unmixed:real}}

\begin{proof}[Proof of Theorem~\ref{unmixed:real}]
Let $\eps > 0$. As in the mixed case, we define:
\begin{eqnarray*}
\nu_{\mathbb R}(U,\eps) &\defeq&
\mathrm{Prob}_{f \in \mathcal F} 
\left[
\mymu(f;U) > \eps^{-1}
\right]
\\
&=&
\mathrm{Prob}_{f \in \mathcal F} 
\left[
\exists p \in U : \mathit{ev}(f;p)=0 \text{\ and \ }
d_{\mathbb P} (f, \Sigma_{p}) < \eps
\right]
\end{eqnarray*}
  where now $U \in \mathbb R^n$.

Let $V(\eps) \defeq \{ (f,p) \in \mathbb \mathcal F_{\mathbb R} \times U: 
\mathit{ev}(f;p)=0 \text{\ and \ }
d_{\mathbb P} (f, \Sigma_{p}) < \eps \}$. We also define $\pi: V(\eps)
\rightarrow \mathbb P(\mathcal F)$ to be the canonical projection 
mapping $F_\R\times U$ to $F_\R$ and set 
$\#_{V(\eps)}(f) \defeq \# \{ p \in U: (f,p) \in V(\eps) \}$.
Then,

\begin{eqnarray*}
\nu_{\mathbb R}(U,\eps)
&=&
\int_{f \in \mathcal F^{\mathbb R}}
\frac{e^{-\sum_i \|f^i\|^2/2}}{\sqrt{2\pi}^{\sum M_i}}
\chi_{\pi(V(\eps))} (f) \ 
d\mathcal F^{\mathbb R} 
\\
&\le&
\int_{f \in \mathcal F^{\mathbb R}}
\frac{e^{-\sum_i \|f^i\|^2/2}}{\sqrt{2\pi}^{\sum M_i}}
\#_{V(\eps)} 
d\mathcal F^{\mathbb R} 
\\
&\le&
\int_{p\in U \subset \mathbb R^n}
\int_{
\substack{
f \in \mathcal F^{\mathbb R}_{p} \\
d_{\mathbb P} (f, \Sigma_{p}) < \eps
}}
\frac{e^{-\sum_i \|f^i\|^2/2}}{\sqrt{2\pi}^{\sum M_i}}
\frac{1}{NJ(f;p)} 
d\mathcal F^{\mathbb R}_{p} 
\ dV_{\toric{n}}
\end{eqnarray*}

\par
As before, we change coordinates in each fiber of
$\mathcal F^{\mathbb R}_{A}$ by
\[
f = f_{\romone} + f_{\romtwo} + f_{\romthree}
\]
with $f^i_{\romone}$ colinear to $v_{A}^T$,
$(f^i_{\romtwo})^T$ in the range of $Dv_{A}$,
and $f^i_{\romthree}$ othogonal to $f^{i}_{\romone}$
and $f^i_{\romtwo}$.
This coordinate system is dependent on $p+q\sqrt{-1}$.
\par
In the new coordinate system, formula~\ref{det:cond:matrix}
splits as follows: 
\begin{eqnarray*}
\det \left(DG_{\scriptscriptstyle (p)} DG_{\scriptscriptstyle (p)}^H\right)^{-1/2}
dV_{\toric{n}} =
\hspace{-97pt} && \\
&=&
\left| \det 
\left[
\begin{matrix}
(f^1_{\romtwo})_1 & \hdots & (f^1_{\romtwo})_n\\ 
\vdots & & \vdots\\
(f^n_{\romtwo})_1 & \hdots & (f^n_{\romtwo})_n\\ 
\end{matrix}
\right]
\right|
\left| \det 
\left[
\begin{matrix}
({Dv_{A}}^{\romtwo})^1_1 & \hdots & ({Dv_{A}}^{\romtwo})^1_n  \\ 
\vdots & & \vdots \\
({Dv_{A}}^{\romtwo})^n_1 & \hdots & ({Dv_{A}}^{\romtwo})^n_n  \\ 
\end{matrix}
\right]
\right|
dV
\\
&=&
\left| \det 
\left[
\begin{matrix}
(f^1_{\romtwo})_1 & \hdots & (f^1_{\romtwo})_n\\ 
\vdots & & \vdots\\
(f^n_{\romtwo})_1 & \hdots & (f^n_{\romtwo})_n\\ 
\end{matrix}
\right]
\right|
\sqrt{\det Dv_A^H Dv_A}
\end{eqnarray*}

The integral $E(U)$ of $\sqrt{\det Dv_A Dv_A^H}$ is the expected number
of real roots on $U$, 
therefore 
\begin{multline*}
\nu_{\mathbb R}(U,\eps)
\le
E(U)
\int_{
\substack{
f_{\romtwo}+f_{\romthree} \in \mathcal F^{\mathbb R}_{p} \\
d_{\mathbb P} (f_{\romtwo}+f_{\romthree}, \Sigma_{p}) < \eps
}}
\frac{e^{-\sum_i \|f_{\romtwo}^i + f_{\romthree}^i\|^2/2}}{\sqrt{2\pi}^{\sum M_i}}
\cdot
\\
\cdot
\left| \det 
\left[
\begin{matrix}
(f^1_{\romtwo})_1 & \hdots & (f^1_{\romtwo})_n\\ 
\vdots & & \vdots\\
(f^n_{\romtwo})_1 & \hdots & (f^n_{\romtwo})_n\\ 
\end{matrix}
\right]
\right|
\ d\mathcal F^{\mathbb R}_{p} \ 
\punct{.}
\end{multline*}

In the new system of coordinates,
$\Sigma_{p}$ is defined by the 
equation:
\[
\det 
\left[
\begin{matrix}
(f^1_{\romtwo})_1 & \hdots & (f^1_{\romtwo})_n\\ 
\vdots & & \vdots\\
(f^n_{\romtwo})_1 & \hdots & (f^n_{\romtwo})_n\\ 
\end{matrix}
\right]
=
0
\punct{.}
\]

Since $\| f_{\romtwo} + f_{\romthree} \| \ge \| f_{\romtwo}\|$,
\[
d_{\mathbb P} (f_{\romtwo}+f_{\romthree}, \Sigma_{p}) < \eps
\Longrightarrow
d_{\mathbb P} (f_{\romtwo}, \Sigma_{p}) < \eps
\punct{.}
\]

This implies:
\begin{multline*}
\nu_{\mathbb R}(U,\eps)
\le
E(U) 
\int_{
\substack{
f_{\romtwo}+f_{\romthree} \in \mathbb \mathcal F^{\mathbb R}_{p} \\
d_{\mathbb P} (f_{\romtwo}, [\det = 0]) < \eps
}}
\frac{e^{-\sum_i \|f_{\romtwo}^i + f_{\romthree}^i\|^2/2}}{\sqrt{2\pi}^{\sum M_i}}
\cdot
\\
\cdot
\left| \det 
\left[
\begin{matrix}
(f^1_{\romtwo})_1 & \hdots & (f^1_{\romtwo})_n\\ 
\vdots & & \vdots\\
(f^n_{\romtwo})_1 & \hdots & (f^n_{\romtwo})_n\\ 
\end{matrix}
\right]
\right|
\ d\mathcal F^{\mathbb R}_{p} \ 
\punct{.}
\end{multline*}

We can integrate the $(\sum M_i - n - 1)$ variables $f_{\romthree}$ to
obtain:
\[
\nu_{\mathbb R}(U,\eps)
=
E(U)
\int_{\substack{
f_{\romtwo} \in \mathbb R^{n^2}\\
d_{\mathbb P} (f_{\romtwo}, [\det = 0]) < \eps
}}
\frac{e^{-\sum_i \|f_{\romtwo}^i\|^2/2}}{\sqrt{2\pi}^{n^2}}
\left| \det f_{\romtwo} \right|^2
\ d\mathbb R^{n^2} 
\punct{.}
\]

This is $E(U)$ times the probability
$\nu(n,\eps)$ for the linear case. 
\end{proof}

\section{Mixed Thoughts about Mixed Manifolds}
\label{finsler}

Let $X: E \rightarrow F$ be a linear operator.
Here, we assume that $E$ has a canonical 
Riemannian structure, and that $F$ has $n$
possibly different Riemannian structures
$\langle \cdot,\cdot\rangle_{\scriptscriptstyle A_i}$.
\par
We would like to interpret the quantities
\[
\max_{\|v\|\le 1} 
\min_j
\| X v\|_{A_i}
\]
and
\[
\max_{\|v\|\le 1} 
\max_j
\| X v\|_{A_i}
\]
in terms of an operator norm on $X$. Let $B$ denote
the unit ball in $E$, and let $B_i$ denote the
unit ball in $(F, \langle\cdot,\cdot\rangle_{\scriptscriptstyle A_i})$. 
Note that 
\par
\[
\max_{\|v\|\le 1} 
\max_j
\| X v\|_{A_i}
\le 1
\
\Longleftrightarrow
\
X(B) \subseteq \bigcap B_i
\punct{,}
\]

while

\[
\max_{\|v\|\le 1} 
\min_j
\| X v\|_{A_i}
\le 1
\
\Longleftrightarrow
\
X(B) \subseteq \bigcup B_i
\punct{.}
\]

This should be compared to:
\[
\min_j
\max_{\|v\|\le 1} 
\| X v\|_{A_i}
\le 1
\
\Longleftrightarrow
\
\exists i: X(B) \subseteq B_i
\punct{.}
\]

One standard way to define norms is by
choosing an arbitrary symmetric convex set,
and equating that set to the unit ball. 
(Such norms are called {\bf Minkowski norms}.) 

There are two immediate obvious choices:

\begin{enumerate}
\item We can use $\bigcap B_i$ as the unit ball.
\item We can use $\cvx \bigcup B_i$ as the unit ball.
\end{enumerate}

In the first case, we can endow \toric{n}
with a $\mathcal C^0$ Finsler structure, while in
the second case we can obtain a $\mathcal C^1$
Finsler structure. 

Using $\cvx \bigcup B_i$ would have the advantage 
of a known probabilistic bound for 
$\mymu > \eps^{-1}$. However, $\bigcap B_i$ seems
to be more convenient for the study of polyhedral
homotopy~\cite{HS}.

Finsler structures are legitimate
ways to endow a non-Riemannian manifold with a 
few familiar concepts. For instance, once we define
a Finsler structure $\finsler{\cdot}_{x}$, the length of a curve $x(t)$,
$t \in [0,1]$ is defined to be $\int_0^1 \finsler{\dot x(t)}_{x(t)} dt$.

A general discussion on Finsler geometry can be found 
in~\cite[Ch.~8]{CHERN-CHEN-LAM}.

\begin{remark}
The proof of Theorem~\ref{condnumber} strongly suggests
that the geometry of mixed manifolds should be determined by a
much more fundamental invariant, a norm in the space
$L(\mathbb C^n,T_{(p,q)}\toric{n})$, which we can take to be 
either side of the following equality: 
\[
\max_v \left( \sqrt{ \sum_{i=1}^n \frac{ |v^i|^2 }{ \| X v \|_{A_i}^2 } 
}\right)^{-1} =
\max_u \left( \sqrt{ \sum_{i=1}^n \frac{ |(X^{-1} u)^i|^2 }{ \| u \|_{A_i}^2 } }\right)^{-1}
\punct{. \qed}
\]
\end{remark}

\begin{remark}
There is a class of polynomial systems that are not unmixed,
but nevertheless can be treated as if they were unmixed.
For instance, in the dense case, the potentials
$g_{A_i}$, $i=1, \cdot, n$ are all multiples of one another,
therefore $\kappa \equiv 1$. The toric variety associated 
to those systems admits therefore a (possibly singular)
Hermitian structure. That structure is non-singular
provided that the $A_i$'s satisfy Delzant's condition
~\cite{DELZANT} (see also Appendix~\ref{momentum} below). 
Roughly speaking, Delzant's condition is an assertion about  
the angle cones of the Minkowski sum of the $\cvx(A_i)$. 
\qed 
\end{remark}

\appendix

\section{Mechanical Interpretation of the Momentum Map}
\label{momentum}

The objective of this Section is to clarify the 
analogy between the geometry of polynomial roots
and Hamiltonian mechanics. The key for that analogy
was the existance of a {\em momentum map} associated
to convex sets. 
\par
In the case the convex set is the support of a polynomial,
that {\em momentum map} is also the momentum map associated
to a certain Lie group action,
namely the natural action of the
$n$-torus on the {\em toric} manifold \toric{n}:

The $n$-torus $\mathbb T^n = \mathbb R^n \pmod {2 \pi \,\mathbb Z^n}$ acts on
$\toric{n}$ by
\[
\rho : (p,q) \mapsto (p,q+\rho) 
\punct{,}
\]
where $\rho \in \mathbb T^n$.

This action preserves the symplectic structure, since it fixes
the $p$-variables and translates the $q$-variables (see Formula~\ref{211}). 
Also,
the Lie algebra of $\mathbb T^n$ is $\mathbb R^n$. An element $\xi$
of $\mathbb R^n$ induces an {\em infinitesimal} action
(i.e. a vector field) $X_{\xi}$ in \toric{n}.

This vector field is the derivation that to any smooth function
$f$ associates:
\[
(X_{\xi})_{(p,q)}(f) = \iota_{\xi} (\omega_A)_{(p,q)} (df)
\defeq
(\omega_A)_{(p,q)} (\xi, df)
\punct{.}
\]

If we write $df = d_p f dp + d_q f dq$, then this formula translates to:

\[
(X_{\xi})_{(p,q)}(f) 
= 
- \xi^T (D^2g_A)_p d_pf
\] 
by Formula~\ref{omega}.

This vector field is Hamiltonian: if $(p(t), q(t))$ is a solution
of the equation
\[
(\dot p(t), \dot q(t)) = (X_{\xi})_{p(t),q(t)}
\]
then we can write
\[
\left\{
\begin{array}{ccl}
\dot p &=& \frac{\partial H_{\xi}}{\partial q}\\
\dot q &=& - \frac{\partial H_{\xi}}{\partial p}\\
\end{array}
\right.
\punct{,}
\]
where $H_{\xi} = \nabla g_A (p) \cdot \xi$.

This construction associates to every $\xi \in \mathbb R^n$, the
Hamiltonian function $H_{\xi} = \nabla g_A (p) \cdot \xi$. The term 
$\nabla g_A(p)$ is a function of $p$, with values in $(\mathbb R^n)^{\vee}$
(the dual of $\mathbb R^n$). In more general Lie group actions,
the momentum map takes values in the dual of the Lie algebra, so
that the pairing $\nabla g_A (p) \cdot \xi$ always makes sense.
A Lie group action with such an expression for the Hamiltonian
is called {\em Hamiltonian} or {\em Strongly Hamiltonian}.

\par

The theorem of Atiyah, Guillemin and Sternberg asserts that
under certain conditions, the image of the momentum map is 
the convex hull of the points $A^{\alpha}$. One of those
conditions requires that the action should take place in a
compact symplectic manifold (as in ~\cite[Th.~1]{ATIYAH82}
or~\cite[Th.~4]{GS}) or (sometimes) a compact K\"ahler 
manifold~\cite[Th~2]{ATIYAH82}.
\par
We may consider a compactification of $\toric{n}$, such 
as the closure of $\exp(\toric{n})$.
Unfortunately,
there are situations where this compactification is not a
symplectic manifold, because the form $\omega_A$ vanishes
in the preimage of one or more $A^{\alpha}$'s (thus differing 
from the assumptions of \cite[p.~127]{ATIYAH83}). 
\par
In~\cite{DELZANT}, a necessary and sufficient condition for 
the compactification of $\toric{n}$ to be a symplectic manifold
is given. Namely, the polytope $\cvx(A)$ should be simple
(i.e., every vertex should be incident to exactly
$n$ edges) and unimodular (the integer points along the 
rays generated by each such $n$-tuple of edges should 
span $\mathbb Z^n$ as a $\Z$-module). If all the polytopes $\cvx(A_1), 
\ldots,\cvx(A_n)$ satisfy Delzant's condition then we can construct a 
corresponding compactifaction 
\toric{n} and apply Atiyah-Guillemin-Sternberg's theorem.
\par
Another possibility is to blow-up the singularities, as explained
in ~\cite{MCDUFF-SALAMON}. If we do so, polytopes $A_1, \cdots, A_n$
will be ``shaved:'' locally, the cone emanating from a vertex 
will be truncated by the intersection with a half-space with 
boundary infinitesimally close to a supporting hyperplane.  
From another point of view, the underlying normal fan~\cite{Ewald}
will be refined. However, the relation between the original 
polynomial system and the new momentum map is not yet clear.

\section{The Coarea Formula}

This is an attempt to give a short proof of the coarea formula,
in a version suitable to the setting of this paper. This means
we take all manifolds and functions smooth and avoid measure
theory as much as possible.

\begin{proposition} \label{coarea2}
\begin{enumerate}
\item
  Let $X$ be a smooth Riemann manifold, of dimension $M$ and 
volume form $|dX|$. 
\item
  Let $Y$ be a smooth Riemann manifold, of
dimension $n$ and volume form $|dY|$. 
\item
  Let $U$ be an open set of $X$, and $F: U \rightarrow Y$
be a smooth map, such that $DF_x$ is
surjective for all $x$ in $U$. 
\item
  Let $\varphi: X \rightarrow \mathbb R^+$ be a 
smooth function with compact support contained in $U$. 
\end{enumerate}
Then for almost all $z \in F(U)$, $V_z \defeq F^{-1}(z)$ is a smooth 
Riemann manifold, and
\[
\int_X \varphi(x) NJ(F;x) |dX|
=
\int_{z \in Y} 
\int_{x \in V_z}
\varphi(x) |dV_z| |dY|
\]
where $|dV_z|$ is the volume element of $V_z$ and
$NJ(F,x) = \sqrt{\det DF_x^H DF_x}$ is the product 
of the singular values of $DF_x$. \qed 
\end{proposition}

 By the implicit function theorem, whenever $V_z$ is non-empty, 
it is a smooth $(N-n)$-dimensional Riemann submanifold
of $X$. By the same reason, 
$V \defeq \{ (z,x): x \in V_z \}$
is also a smooth manifold. 

Let $\eta$ be the following $N$--form restricted to $V$:
\[
\eta = dY \wedge dV_z
\punct{.}
\]
\par
This is {\bf not} the volume form of $V$. 
The proof of Proposition~\ref{coarea2} is divided into two
steps:

\begin{lemma}\label{coarea21}
\[
\int_V \varphi(x) |\eta| = \int_X \varphi(x) NJ(F;x) |dX|
\]
\punct{.}
\end{lemma}

\begin{lemma}\label{coarea22}
\[
\int_V \varphi(x) |\eta| = \int_{z \in Y} \int_{x \in V_z}
\varphi(x) |dV_z| |dY|
\punct{.}
\]
\end{lemma}

\begin{proof}[Proof of Lemma~\ref{coarea21}]
We parametrize:
\[
\begin{array}{lrcl}
\psi: & X & \rightarrow & V \\
&x & \mapsto & (F(x),x)
\end{array}
\punct{.}
\]

Then,
\[
\int_V \varphi (x) |\eta| = \int_X (\varphi \circ \psi) (x) |\psi ^* \eta|
\punct{.}
\]

We can choose an orthonormal basis $u_1, \cdots, u_M$ of $T_xX$ such that
$u_{n+1}, \cdots, u_M \in \ker DF_x$. Then,
\[
D\psi (u_i) = 
\left\{
\begin{array}{ll}
 (DF_x u_i, u_i) & i=1, \cdots, n \\
 (0,u_i)         & i=n+1, \cdots, M
\end{array}
\right.
\punct{.}
\]
Thus,
\begin{eqnarray*}
| \psi^* \eta (u_1, \cdots, u_M) | &=&
| \eta (D \psi u_1, \cdots, D\psi u_M) | \\
&=&
|dY( DF_x u_1, \cdots, DF_x u_n )| 
\ |dV_z(u_{n+1}, \cdots, u_M)|
\\
&=& 
|\det DF_x |_{\ker DF_x^{\perp}}| \\
&=&
NJ(F,x) 
\end{eqnarray*}
and hence
\[
\int_V \varphi(x) |\eta| = \int_X \varphi(x) NJ(F;x) |dX|
\punct{.}
\]
\end{proof}

\begin{proof}[Proof of Lemma~\ref{coarea22}]
We will prove this Lemma locally, and this 
implies the full Lemma through a standard
argument (partitions of unity in a compact
neighborhood of the support of $\varphi$). 
\medskip
\par
Let $x_0,z_0$ be fixed. A small enough neighborhood of
$(x_0,z_0) \subset V_{z_0}$ admits a fibration
over $V_{z_0}$ 
by planes orthogonal to $\ker DF_{x_0}$.
\par
We parametrize:
\[
\begin{array}{llcl}
\theta: & Y \times V_{z_0} & \rightarrow & V \\
& (z,x) & \mapsto & (z, \rho(x,z))
\end{array}
\punct{,}
\]

where $\rho(x,z)$ is the solution of $F(\rho)=z$ in the fiber
passing through $(z_0,x)$.
Remark that $\theta^* dY = dY$, and $\theta^* dV_z = \rho^* DV_z$.
Therefore,
\[
\theta^* (dY \wedge dV_z) = dY \wedge (\rho^* dV_{z})
\punct{.}
\]

Also, if one fixes $z$, then $\rho$ is a parametrization
$V_{z_0} \rightarrow V_z$. We have:
\begin{eqnarray*}
\int_V \varphi(x) |\eta| 
&=& 
\int_{Y \times V_{z_0}} \varphi(\rho(x,z)) |\theta^* \eta|
\\
&=& 
\int_{z \in Y} 
\left(
\int_{x \in V_{z_0}} \varphi(\rho(x,z) |\rho^* dV_{z}|
\right) |dY|
\\
&=&
\int_{z \in Y} 
\left(
\int_{x \in V_{z}} \varphi(x) |dV_{z}|
\right) |dY|
\end{eqnarray*}
\end{proof}

The proposition below is essentially Theorem~3 p.~240 of~\cite{BCSS}.
However, we do not require our manifolds to be compact. We assume
all maps and manifolds are smooth, so that we can apply 
proposition~\ref{coarea2}.

\begin{proposition}\label{coarea1}\mbox{}\\
\vspace{-.5cm} 
\begin{enumerate}
\item \label{ca1} Let $X$ be a smooth $M$-dimensional manifold with volume element $|dX|$.
\item \label{ca2} Let $Y$ be a smooth $n$-dimensional manifold with volume element $|dY|$.
\item \label{ca3} Let $V$ be a smooth $M$-dimensional submanifold of $X \times Y$, and
     let $\pi_1: V \rightarrow X$ and $\pi_2: V \rightarrow Y$ be the 
canonical projections from $X\times Y$ to its factors. 
\item \label{ca4} Let $\Sigma'$ be the set of critical points of $\pi_1$,
     we assume that $\Sigma'$ has measure zero and that $\Sigma'$ is
     a manifold.
\item \label{ca5} We assume that $\pi_2$ is regular (all points in $\pi_2(V)$ are regular
     values).
\item \label{ca6} For any open set $U \subset V$, for any $x \in X$, 
     we write: $\#_U(x) \defeq \# \{ \pi_1^{-1}(x) \cap U \}$. We
     assume that $\int_{x \in X} \#_V (x) |dX|$ is finite. 
\end{enumerate}
\par
     Then, for any open set $U \subset V$,
\[
\int_{x \in \pi_1(U)} \#_U (x) |dX| 
=
\int_{z \in Y}
\int_{\substack{x \in V_z\\ (x,z) \in U}} \frac{1}{\sqrt{ \det DG_x DG_x^H }} |dV_z| |dY|  
\]
  where $G$ is the implicit function for $(\hat x,G(\hat x)) \in V$ in a 
neighborhood of $(x,z) \in V \setminus \Sigma'$. \qed 
\end{proposition}

\begin{proof}
  Every $(x,z) \in U \setminus \Sigma'$ admits an open neighborhood such that
$\pi_1$ restricted to that neighborhood is a diffeomorphism. This defines
an open covering of $U \setminus \Sigma'$. Since $U \setminus \Sigma'$ is
locally compact, we can take a countable subcovering and define a partition of
unity $(\varphi_{\lambda})_{\lambda \in \Lambda}$ subordinated to that 
subcovering.
\par
Also, if we fix a value of $z$, then $(\varphi_{\lambda})_{\lambda \in \Lambda}$
becomes a partition of unity for $\pi_1(\pi_1^{-1} (V_z) \cap U)$. Therefore,
\begin{eqnarray*}
\int_{x \in \pi_1(U)} \#_U (x) |dX| 
&=&
\sum_{\lambda \in \Lambda}
\int_{x,z \in \supp \varphi_{\lambda}} \varphi_{\lambda}(x,z) |dX| 
\\
&=&
\sum_{\lambda \in \Lambda}
\int_{z \in Y}
\int_{x,z \in \supp \varphi_{\lambda}} \frac{\varphi_{\lambda}(x,z)}{NJ(G,x)} |dX| 
\\
&=&
\int_{z \in Y}
\sum_{\lambda \in \Lambda}
\int_{x,z \in \supp \varphi_{\lambda}} \frac{\varphi_{\lambda}(x,z)}{NJ(G,x)} |dX| 
\\
&=&
\int_{z \in Y}
\int_{x \in V_z} \frac{1}{NJ(G,x)} |dX| 
\end{eqnarray*}

where the second equality uses Proposition~\ref{coarea2} with $\varphi = \varphi_{\lambda} / NJ$.
Since $NJ = \sqrt{ \det DG_x DG_x^H }$, we are done.
\end{proof}

\bibliographystyle{amsalpha}
{
\footnotesize
\bibliography{poly}
}

\end{document}